\documentclass[11pt,twoside]{article}
\usepackage[hmargin=0.8in,vmargin=0.9in]{geometry}
\usepackage{fancyhdr}
\usepackage{graphicx}
\usepackage{subfigure}
\usepackage{amssymb}
\usepackage{amsmath}
\usepackage{amsfonts}
\usepackage{theorem}
\usepackage{mathrsfs}
\usepackage{mathtools}
\usepackage{bm}
\usepackage{color}
\usepackage{setspace}
\usepackage{exscale}
\usepackage{relsize}
\usepackage{tikz}
\usetikzlibrary{positioning}
\usepackage{rotating}
\usepackage{enumerate}
\usepackage{hyperref}
\usepackage{authblk}
\usepackage{epstopdf}
\usepackage{float}
\usepackage{cite}

\usepackage{cleveref}

\usepackage{nicefrac}
\usepackage{setspace}

\usepackage{booktabs,multirow} 
\usepackage{array} 
\usepackage{paralist} 
\usepackage{verbatim} 
\usepackage{subfigure} 

\theoremstyle{plain}			
\newtheorem{thm}{Theorem}[section]

\newtheorem{rmk}[thm]{Remark}

{\theorembodyfont{\rmfamily}}

\newenvironment{DA}{{\flushleft \bf Declarations:}}{}

\usepackage{tikz}
\usetikzlibrary{positioning}
\usepackage{xcolor}

\allowdisplaybreaks[1]

\numberwithin{equation}{section}
\numberwithin{figure}{section}
\numberwithin{table}{section}

\newcommand\eref[1]{(\ref{#1})}

\newcommand*\xbar[1]{%
  \hbox{%
    \vbox{%
      \hrule height 0.5pt 
      \kern0.4ex
      \hbox{%
        \kern-0.05em
        \ensuremath{#1}%
        \kern-0.00em
      }%
    }%
  }%
}

\newcommand{\bmR}{\bm{\mathcal{R}}}
\newcommand{\mF}{\bm{F}}

\newcommand{\mG}{\bm{G}}

\newcommand{\mU}{\bm{U}}

\newcommand{\dt}{\Delta t}
\newcommand{\dx}{\Delta x}
\newcommand{\dy}{\Delta y}

\newcommand{\hf}{{\frac{1}{2}}}

\newcommand{\jph}{{j+\frac{1}{2}}}
\newcommand{\jmh}{{j-\frac{1}{2}}}
\newcommand{\kph}{{k+\frac{1}{2}}}
\newcommand{\kmh}{{k-\frac{1}{2}}}

\def\softd{{\leavevmode\setbox1=\hbox{d}%
          \hbox to 1.05\wd1{d\kern-0.4ex{\char039}\hss}}}

\title{Novel Adaptive Methods for Hyperbolic Conservation Laws Based on New Quasi-Linear Seventh- and Ninth-Order Schemes}
\author{Shaoshuai Chu\thanks{Department of Mathematics, RWTH Aachen University, 52056 Aachen, Germany; {\tt chu@igpm.rwth-aachen.de}},
Pingyao Feng\thanks{Department of Mathematics, North Carolina State University, 27695 Raleigh, USA; {\tt pfeng3@ncsu.edu}}, Vadim A.
Kolotilov\thanks{Lavrentyev Institute of Hydrodynamics Siberian Branch of the Russian Academy of Sciences, Novosibirsk, 630090 Russia and
Keldysh Institute of Applied Mathematics of the Russian Academy of Sciences, 125047 Moscow, Russia; {\tt kolotilov1992@gmail.com}},\\
Alexander Kurganov\thanks{Department of Mathematics and Shenzhen International Center for Mathematics, Southern University of Science and
Technology, Shenzhen, 518055, China; {\tt alexander@sustech.edu.cn}}, and Vladimir V. Ostapenko\thanks{Lavrentyev Institute of Hydrodynamics
Siberian Branch of the Russian Academy of Sciences, Novosibirsk, 630090 Russia and Keldysh Institute of Applied Mathematics of the Russian
Academy of Sciences, 125047 Moscow, Russia; {\tt ostigil@mail.ru}}}

\begin{document} 

\date{}
\maketitle

\begin{abstract}
We develop new adaptive numerical schemes for one- and two-dimensional  hyperbolic systems of conservation laws. The methodology relies on
the use of a smoothness indicator to automatically partition the computational domain into smooth and nonsmooth (``rough``) regions. We then
follow the scheme adaption strategy recently introduced in [{\sc S. Chu, P. Feng, V. A. Kolotilov, A. Kurganov, and V. V. Ostapenko},
Commun. Comput. Phys., accepted], but instead of the quasi-linear (QL) fifth-order finite-difference scheme used there, we employ the new QL
seventh- and ninth-order schemes in the smooth regions. A series of numerical experiments for the Euler equations of gas dynamics
demonstrates that the new adaptive schemes contain a smaller amount of numerical dissipation and achieve higher resolution compared with
their counterpart that uses the QL fifth-order scheme in the smooth areas. 
\end{abstract}



\section{Introduction}
In this paper, we consider hyperbolic systems of conservation laws, which in the one- (1-D) and two-dimensional (2-D) cases, read as
\begin{equation}
\mU_t+\mF(\mU)_x=\bm0,
\label{1.1}
\end{equation}
and
\begin{equation}
\mU_t+\mF(\mU)_x+\mG(\mU)_y=\bm0,
\label{1.2}
\end{equation}
respectively. Here, $x$ and $y$ are spatial variables, $t$ is the time, $\mU\in\mathbb R^d$ is a vector of unknown functions, and
$\mF:\mathbb R^d\to\mathbb R^d$ and $\mG:\mathbb R^d\to\mathbb R^d$ are fluxes.

It is well-known that solutions of \eref{1.1} and \eref{1.2} may develop complex structures consisting of shocks, rarefactions, and contact
discontinuities. This complexity makes it difficult to simultaneously achieve (i) sharp, non-oscillatory resolution of discontinuities and
(ii) high-order accuracy in smooth regions. For a variety of existing numerical methods for \eref{1.1} and \eref{1.2}, we refer the reader
to the monographs and review papers \cite{GR2,Tor,Leveque02,Hesthaven18,BAF,Shu20} and references therein. Most of those methods rely on
computationally expensive nonlinear limiters, which are, in fact, not needed to be applied in smooth regions. 

In this paper, we focus on highly efficient scheme adaption methods, which are based on automatic detection of smooth and nonsmooth
(``rough'') parts of the computed solutions and on applying different numerical methods in the corresponding parts of the computational
domain. In particular, we follow our recent paper \cite{CFKKO2025}, where we have applied a smoothness indicator (SI) introduced in
\cite{RL87} to dynamically partition the computational domain into three parts: smooth, nonsmooth, and contact areas within the nonsmooth
region. For the Euler equations of gas dynamics, this can be done by applying the SI to the density and pressure variables separately as
pressure remains continuous at the contact surfaces, where the density jumps, while both the density and pressure jump across shock waves.
After identifying the locations of contact surfaces, we apply the second-order low-dissipation central-upwind (LDCU) schemes introduced in
\cite{CKX_24} with an overcompressive limiter from \cite{Lie03} in their vicinities only, while utilizing the dissipative Minmod2 limiter in
the rest of the ``rough'' areas. Both the limiters are applied to the local characteristic variables using the local characteristic
decomposition, which is often used in the context of high-order schemes \cite{Qiu02}, but can also be implemented to suppress spurious
oscillations, which might be otherwise developed by second-order schemes; see, e.g., \cite{CKX24,CKM2025}. Finally, in the smooth regions,
we have used the quasi-linear (QL) fifth-order scheme from \cite{KOK25}, which is highly accurate and computationally inexpensive.  

The main goal of the present paper is to further increase the resolution of the smooth parts of the computed solution. To this end, we
replace the QL fifth-order schemes, which were used in \cite{CFKKO2025}, with the recently developed QL seventh- and ninth-order schemes
from \cite{KOK79}, and incorporate them into the same adaptive strategy. The new QL higher-order schemes rely on a high-order
artificial-viscosity term whose strength is controlled by parameters that can be selected to optimize their linear stability properties. To
determine the optimal values of those parameters, we conduct a von Neumann stability analysis for the 1-D QL seventh- and ninth-order
schemes. In the 2-D extension, we follow a dimension-by-dimension construction and retain the same stability-motivated viscosity parameters.
The resulting adaptive schemes preserve the robustness and sharp resolution of the discontinuities, while providing a higher-fidelity
representation of smooth structures. The improvements are illustrated on a series of 1-D and 2-D benchmarks for the Euler equations of gas
dynamics, where the seventh- and ninth-order schemes are compared against the fifth-order one within the same adaptive framework. 

The rest of the paper is organized as follows. In \S\ref{sec2}, we present the 1-D $r^{\rm th}$-order QL schemes, where $r\ge3$ is an odd
number. For $r=7$ and $9$, we study linear stability of these schemes, verify their orders of accuracy, and introduce the 1-D adaptive
algorithms. In \S\ref{sec3}, we extend the $r^{\rm th}$-order QL schemes to the 2-D case, and for $r=7$ and $9$, we introduce the 2-D
adaptive algorithms. In \S\ref{sec4}, we report numerical experiments demonstrating the performance of the proposed schemes. 

\section{One-Dimensional Schemes}\label{sec2}
In this section, we consider the 1-D hyperbolic system \eref{1.1} and describe the adaption algorithm, which is similar to the one recently
introduced in \cite{CFKKO2025}. We assume that the computational domain is covered with the uniform cells $I_j:=\big[x_\jmh,x_\jph\big]$
with $x_\jph-x_\jmh\equiv\dx$ centered at $x_j=\big(x_\jmh+x_\jph\big)/2$.

\subsection{Quasi-Linear $r^{\rm th}$-Order Schemes}\label{sec21}
We first briefly review the QL $r^{\rm th}$-order schemes (with $r=2\ell-1$, $\ell=2,3,\ldots$) recently proposed in \cite{KOK79}. These
schemes are, in fact, $r^{\rm th}$ order in space and third order in time and they can be written in the following three stages:
\begin{equation*}
\mU^{\rm I}_j=\mU^n_j-\dt^n\bmR^n_j,\quad
\mU^{\rm II}_j=\frac{3}{4}\,\mU_j^n+\frac{1}{4}\left[\mU^{\rm I}_j-\dt^n\bmR^{\rm I}_j\right],\quad
\mU_j^{n+1}=\frac{1}{3}\,\mU_j^n+\frac{2}{3}\left[\mU^{\rm II}_j-\dt^n\bmR^{\rm II}_j\right].
\end{equation*}
Here, $\dt^n$ is a time step, $\mU^n_j\approx\mU(x_j,t^n)$, $t^n$ is the current time level, $t^{n+1}=t^n+\dt^n$, and $\bmR^n_j$,
$\bmR^{\rm I}_j$, and $\bmR^{\rm II}_j$ are the following discretizations of $\mF(\mU)_x$:
\begin{equation}
\bmR^n_j=(\mF_x)^n_j,\quad\bmR^{\rm I}_j=(\mF_x)^{\rm I}_j,\quad
\bmR^{\rm II}_j=(\mF_x)^{\rm II}_j+(-1)^{\ell}\,\frac{3(\dx)^{r+1}}{2\dt^n}B(D^{r+1}_x\mU)^n_j.
\label{2.1a}
\end{equation}
Here, $(\mF_x)^n_j$ is a $(r+1)^{\rm st}$-order central-difference approximation of $\mF(\mU(x_j,t^n))_x$, which is given by 
\begin{equation}
(\mF_x)^n_j=\frac{1}{\dx}\sum_{i=1}^\ell a_i\big(\mF(\mU^n_{j+i})-\mF(\mU^n_{j-i})\big)  
\label{2.1b}
\end{equation}
with the coefficients $a_i$ obtained by solving the following linear algebraic system:
\begin{equation}
\sum_{i=1}^\ell a_ii=\hf,\qquad\sum_{i=1}^\ell a_ii^{2s+1}=0,\quad s=1,\dots,\ell-1.
\label{2.1aa}
\end{equation}
The terms $(\mF_x)^{\rm I}_j$ and $(\mF_x)^{\rm II}_j$ in \eref{2.1a} are determined similarly with the superscript $n$ in \eref{2.1b}
replaced with ${\rm I}$ and ${\rm II}$, respectively.
 
Finally, $(D^{r+1}_x\mU)^n_j$ is the following second-order central-difference approximation of
$\frac{\partial^{r+1}}{\partial x^{r+1}}\mU(x_j,t^n)$:
$$
(D^{r+1}_x\mU)^n_j=\frac{1}{(\dx)^{r+1}}\sum_{i=-\ell}^\ell(-1)^{i+\ell}\binom{2\ell}{i+\ell}\mU^n_{j+i}, 
$$ 
which is a part of the term that represents high-order artificial viscosity used to stabilize the QL schemes. The amount of this artificial
viscosity can be controlled by the positive parameter $B$ in \eref{2.1a}, which can be determined based on the linear stability analysis,
which we conduct for the seventh- and ninth-order schemes in \S\ref{sec2b} and \S\ref{sec2c}, respectively. 
 
It is easy to show that the above QL $r^{\rm th}$-order schemes can be written in the following flux form:
\begin{equation}
\begin{aligned}
\mU^{\rm I}_j&=\mU^n_j-\lambda^n\left(\bm{{\cal F}}_\jph-\bm{{\cal F}}_\jmh\right),\\
\mU^{\rm II}_j&=\frac{3}{4}\,\mU_j^n+\frac{1}{4}\left[\mU^{\rm I}_j-\lambda^n\left(\bm{{\cal F}}^{\rm I}_\jph-\bm{{\cal F}}^{\rm I}_\jmh
\right)\right],\\[0.5ex]
\mU_j^{n+1}&=\frac{1}{3}\,\mU_j^n+\frac{2}{3}\left[\mU^{\rm II}_j-\lambda^n\left(\bm{{\cal F}}^{\rm II}_\jph-\bm{{\cal F}}^{\rm II}_\jmh
\right)\right],
\end{aligned}
\label{2.1bb}
\end{equation}
where $\lambda^n:={\dt^n}/{\dx}$ and the numerical fluxes are defined by:
\begin{equation}
\bm{{\cal F}}_\jph=\bm{{\cal L}}_\jph[\bm U^n],\quad\bm{{\cal F}}^{\rm I}_\jph=\bm{{\cal L}}_\jph\big[\bm U^{\rm I}\big],\quad
\bm{{\cal F}}^{\rm II}_\jph=\bm{{\cal L}}_\jph\big[\bm U^{\rm II}\big]-\bm\omega_\jph^n,
\label{2.1cc}
\end{equation}
where
\begin{equation*}
\bm{{\cal L}}_\jph[\bm U]:=\sum_{i=1-\ell}^\ell b_i\mF(\mU_{j+i})\quad{\rm and}\quad\bm\omega_\jph^n:=\frac{3B}{2\lambda^n}
\sum_{i=1-\ell}^\ell(-1)^{i+\ell}\binom{2\ell-1}{i+\ell-1}\mF(\mU_{j+i}). 
\end{equation*}
Here, the coefficients
\begin{equation}
b_i=b_{1-i}=\sum_{m=i}^\ell a_m,\quad i=1,\dots,\ell,
\label{2.1dd}
\end{equation}
with $a_m$ defined by \eref{2.1aa}.

\subsection{Von Neumann Stability Analysis of the Quasi-Linear Seventh-Order Schemes}\label{sec2b}
\begin{figure}[ht!]
\centerline{\includegraphics[trim=0.1cm 0.1cm 0.1cm 0.1cm, clip, width=4.5cm]{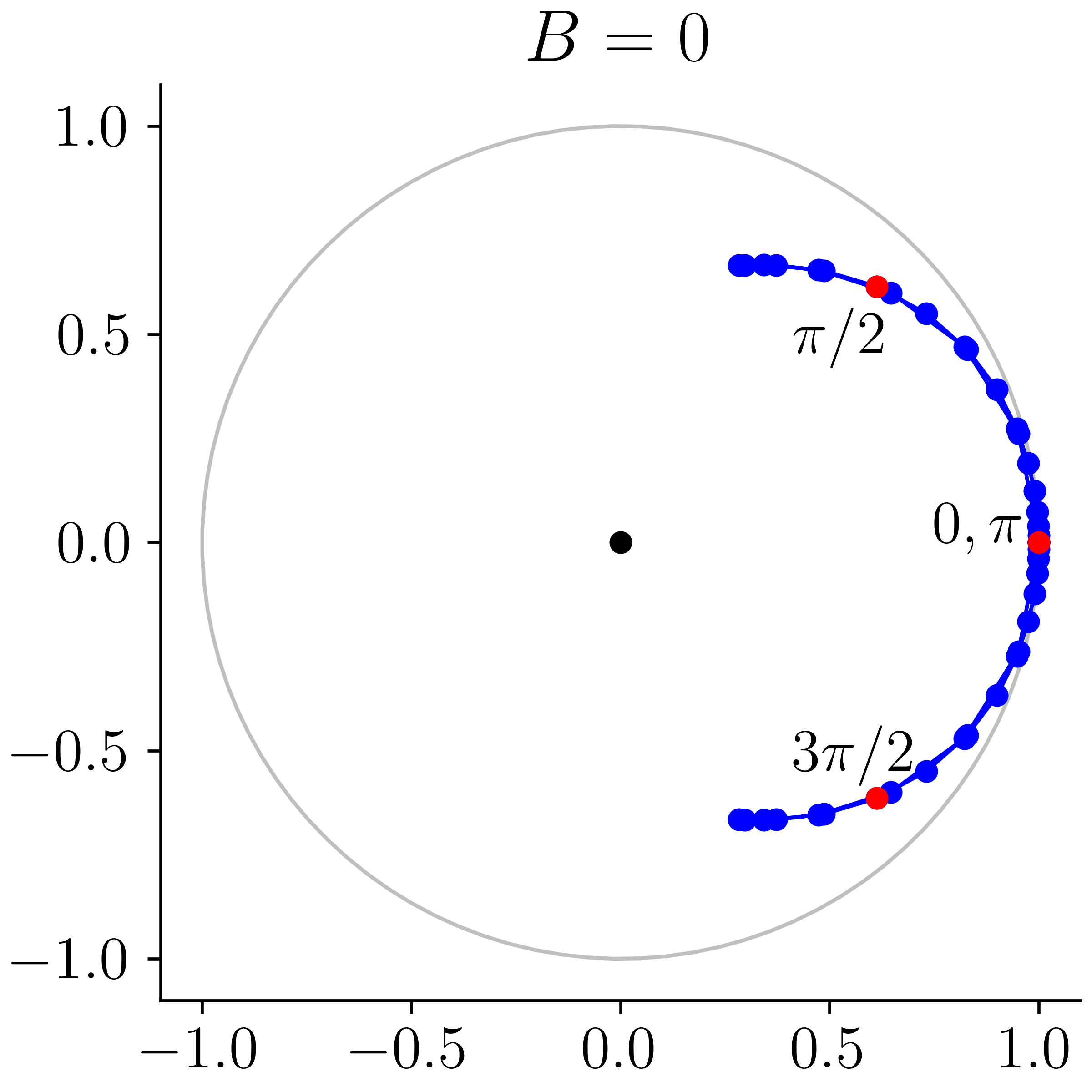}\hspace{1cm}
            \includegraphics[trim=0.1cm 0.1cm 0.1cm 0.1cm, clip, width=4.5cm]{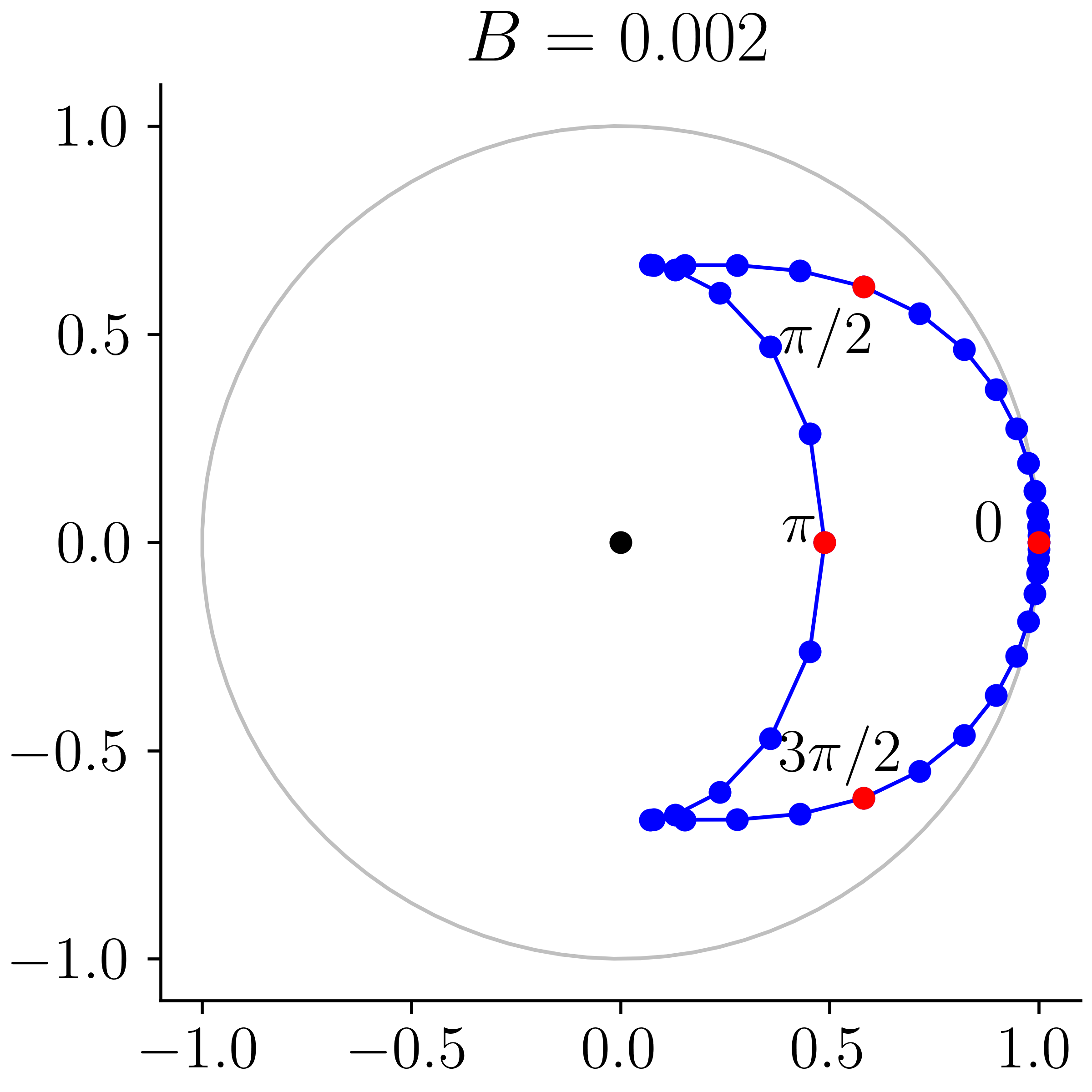}}
\vskip15pt 
\centerline{\includegraphics[trim=0.1cm 0.1cm 0.1cm 0.1cm, clip, width=4.5cm]{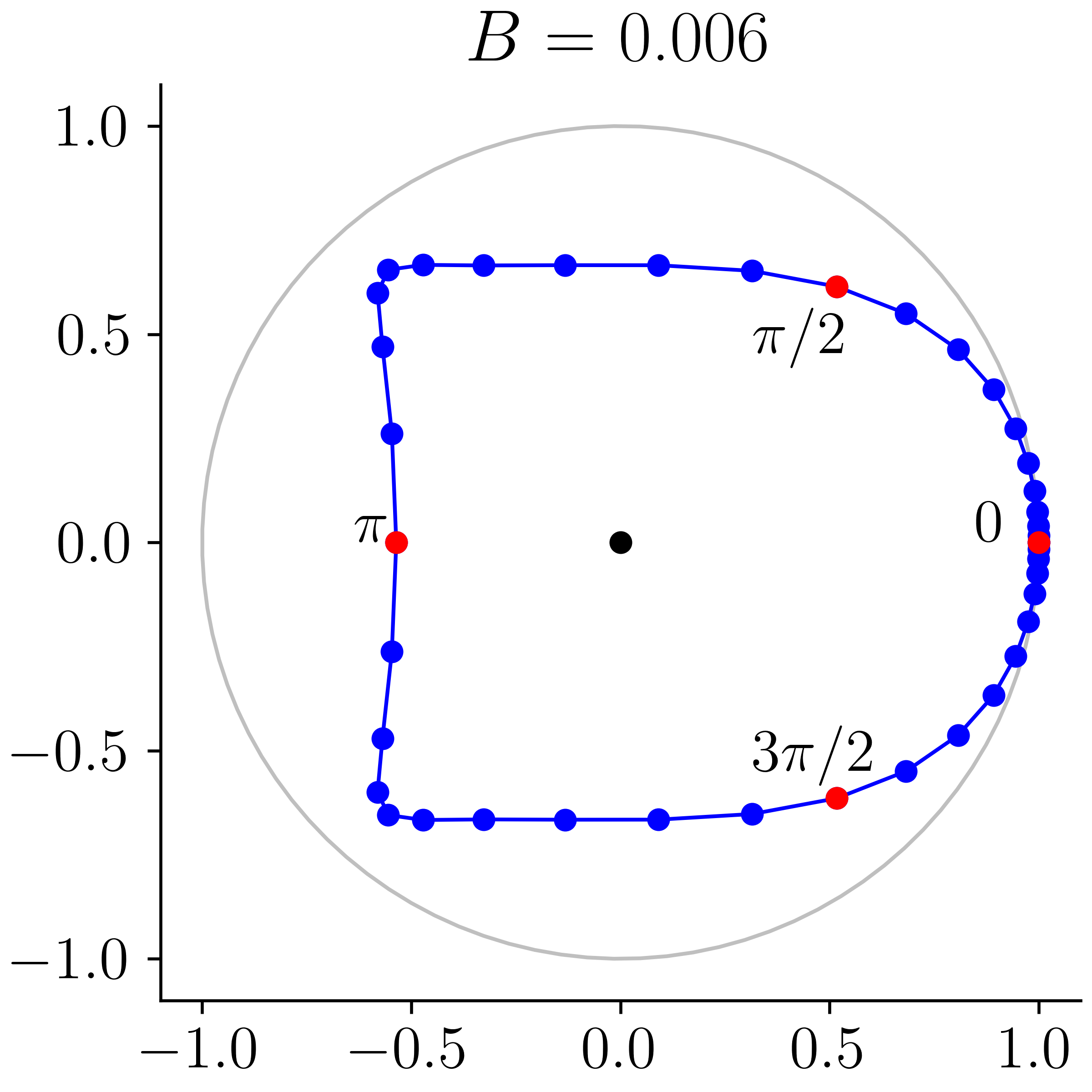}\hspace{1cm}
            \includegraphics[trim=0.1cm 0.1cm 0.1cm 0.1cm, clip, width=4.5cm]{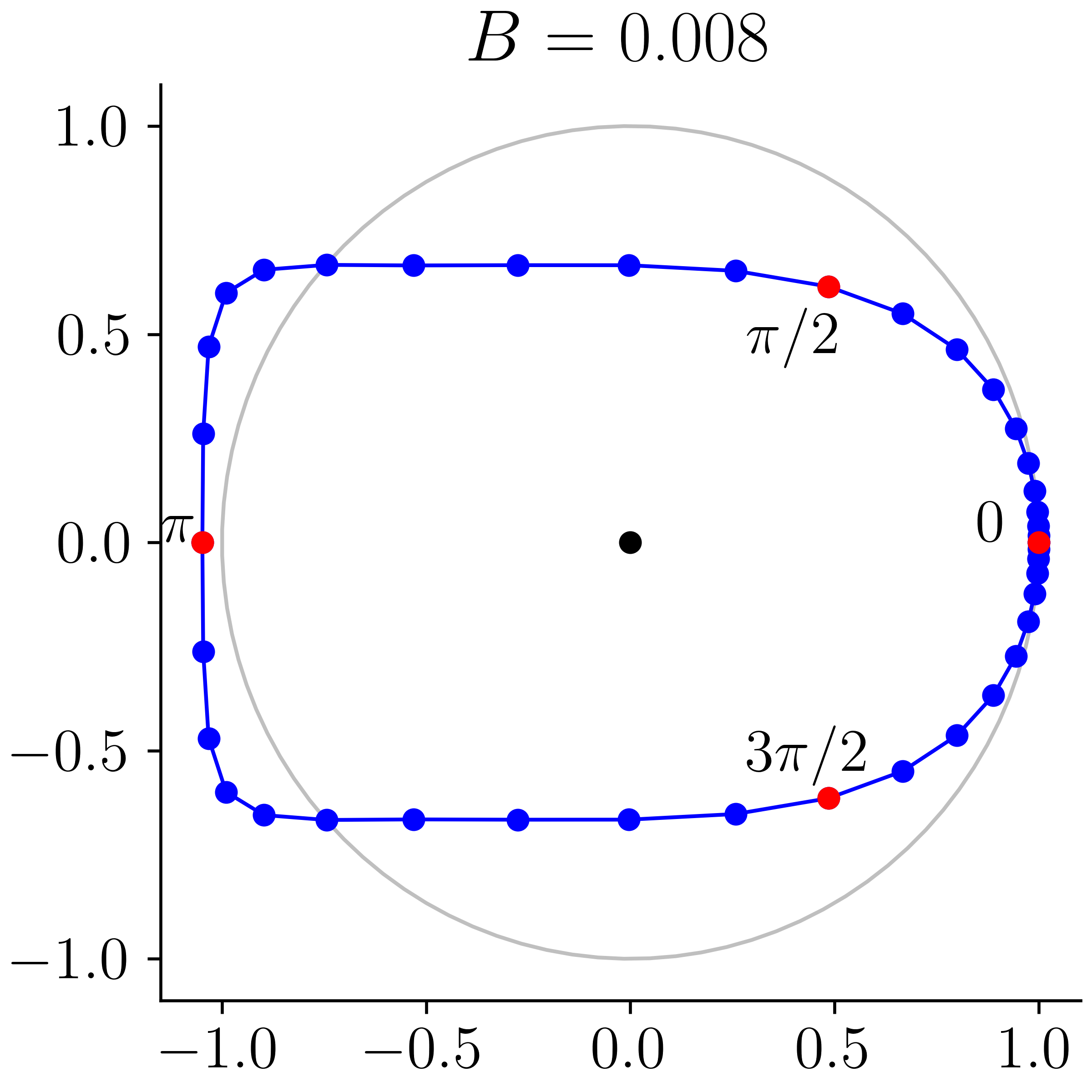}}            
\caption{\sf QL seventh-order scheme: The values $g=g(0.5,B,\alpha)$ relative to the unit circle in the complex plane for
$\alpha=\alpha_\ell$ given by \eref{21aa} and different $B$.\label{fig3}}
\end{figure}

We consider the QL seventh-order schemes, which reads as \eref{2.1bb}--\eref{2.1cc} with
\begin{equation}
\begin{aligned}
&\bm{{\cal L}}_\jph[\bm U]=\frac{1}{840}\big[-3\mF(\mU_{j+4})+29\mF(\mU_{j+3})-139\mF(\mU_{j+2})+533\mF(\mU_{j+1})\\
&\hspace{2.8cm}+533\mF(\mU_j)-139\mF(\mU_{j-1})+29\mF(\mU_{j-2})-3\mF(\mU_{j-3})\big],\\
&\bm\omega_\jph^n=
\frac{3}{512\lambda^n}\big[-\mU_{j+4}^n+7\mU_{j+3}^n-21\mU_{j+2}^n+35\mU_{j+1}^n-35\mU_j^n+21\mU_{j-1}^n-7\mU_{j-2}^n+\mU_{j-3}^n\big],
\end{aligned}
\label{2.1bb1}
\end{equation}
and apply them to the linear advection equation
\begin{equation}
\psi_t+\psi_x=0,
\label{Neum1}
\end{equation}
for which the speeds of propagation do not change and hence $\lambda^n\equiv\lambda$ as the size of time steps is independent of $n$,
namely, $\dt^n\equiv\dt$. We then substitute the Fourier mode
\begin{equation}
\psi_j^n=g^ne^{ij\alpha},\quad i=\sqrt{-1},
\label{Neum2}
\end{equation}
into the scheme \eref{2.1bb}--\eref{2.1cc}, \eref{2.1bb1} and obtain the amplification factor
\begin{equation}
g=g(\lambda,B,\alpha)=\frac{1}{3}+\frac{2}{3}\left(1+\lambda\delta(\alpha)\right)
\left(1+\hf\lambda\delta(\alpha)+\frac{1}{4}(\lambda)^2\delta(\alpha)^2\right)+B\Theta(\alpha),
\label{Neum6}
\end{equation}
where 
\begin{equation}
\begin{aligned}
\delta(\alpha)&=\frac{1}{280}\big(e^{4i\alpha}-e^{-4i\alpha}\big)-\frac{4}{105}\big(e^{3i\alpha}-e^{-3i\alpha}\big)
+\frac{1}{5}\big(e^{2i\alpha}-e^{-2i\alpha}\big)-\frac{4}{5}\big(e^{i\alpha}-e^{-i\alpha}\big)\\
&=i\left(\frac{1}{140}\sin(4\alpha)-\frac{8}{105}\sin(3\alpha)+\frac{2}{5}\sin(2\alpha)-\frac{8}{5}\sin\alpha\right),\\
\Theta(\alpha)&=-(e^{4i\alpha}+e^{-4i\alpha})+8(e^{3i\alpha}+e^{-3i\alpha})-28(e^{2i\alpha}+e^{-2i\alpha})+56(e^{i\alpha}+e^{-i\alpha})-70\\
&=-2\cos(4\alpha)+16\cos(3\alpha)-56\cos(2\alpha)+112\cos(\alpha)-70.
\end{aligned}
\label{Neum7} 
\end{equation}
The QL seventh-order scheme is then stable provided
\begin{equation*}
|g(\lambda,B,\alpha)|\le1,\quad\forall\alpha\in[0,2\pi].
\end{equation*}
First, we note that if $B=0$ then $g(\lambda,B,\pi)=1$, and the resulting scheme will produce substantial oscillations when it is applied to
an initial value problem with discontinuous initial data. We therefore need to use $B>0$. To select a particular value of $B$, we make the
following considerations. In Figure \ref{fig3}, we plot the values of $g$ for $\lambda=0.5$ (this is an upper bound for the CFL number for
the LDCU schemes, which will be used in our adaptive methods) for different values of
\begin{equation}
\alpha=\alpha_\ell=\frac{\pi\ell}{20},\quad\ell=0,\dots,40,
\label{21aa}
\end{equation}
and $B=0$, $0.002$, $0.006$, and $0.008$. As one can see, when $B=0.008$, the scheme is unstable, while for $B=0.002$ and $0.006$, it is
stable. To determine the optimal value of $B$, we follow \cite{KOK25,Ostapenko2000} and choose $B=B_*$, for which 
\begin{equation}
g(\lambda,B_*,\pi)=0.
\label{Neum11}
\end{equation}
The use of $B_*$ will enforce the strongest damping of grid-scale oscillations, which correspond to the highest-frequency mode $\alpha=\pi$.
Substituting \eref{Neum6}--\eref{Neum7} into \eref{Neum11} results in $B_*=1/256$. 

\begin{figure}[ht!]
\centerline{\includegraphics[trim=0.1cm 0.1cm 0.1cm 0.1cm, clip, width=4.5cm]{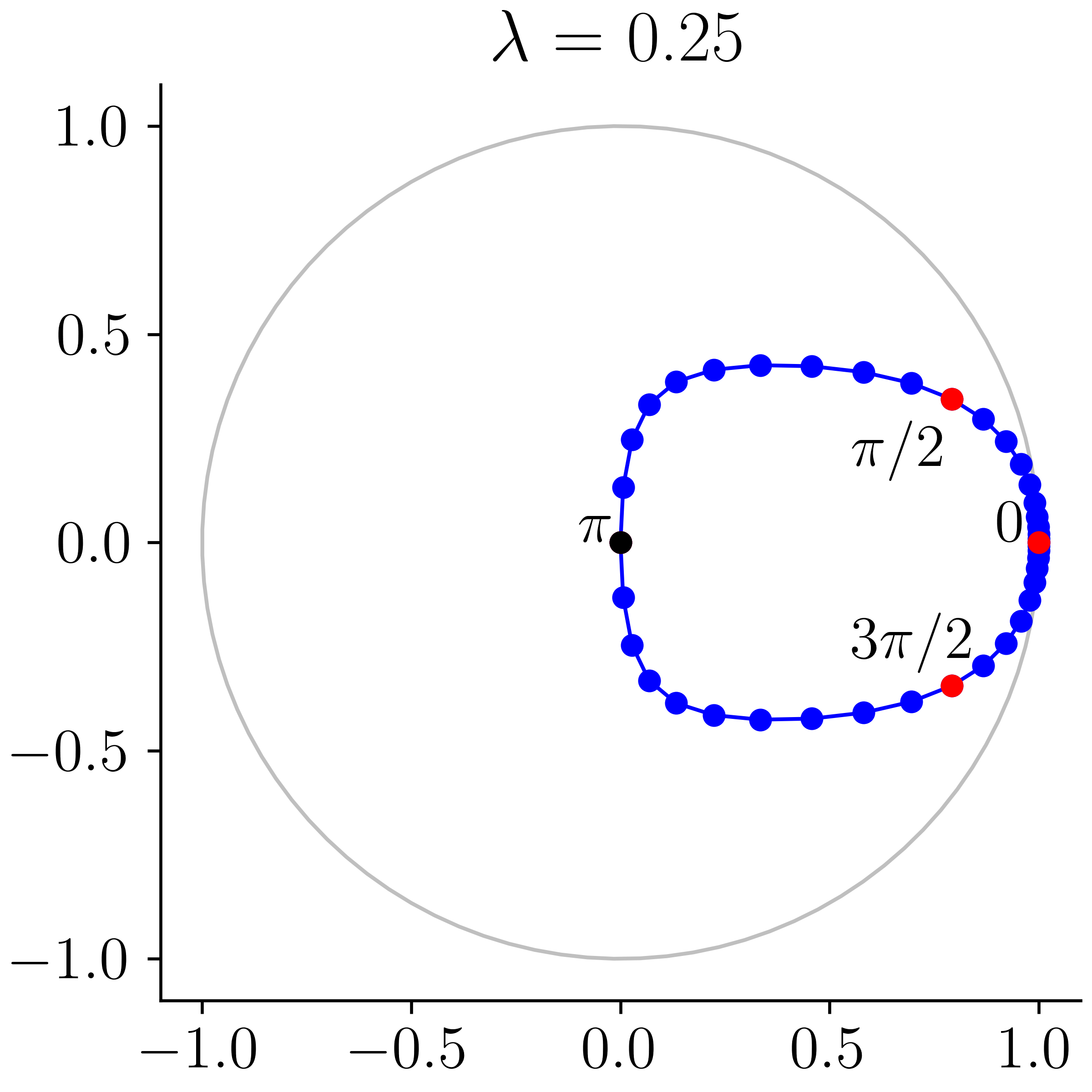}\hspace{1cm}
            \includegraphics[trim=0.1cm 0.1cm 0.1cm 0.1cm, clip, width=4.5cm]{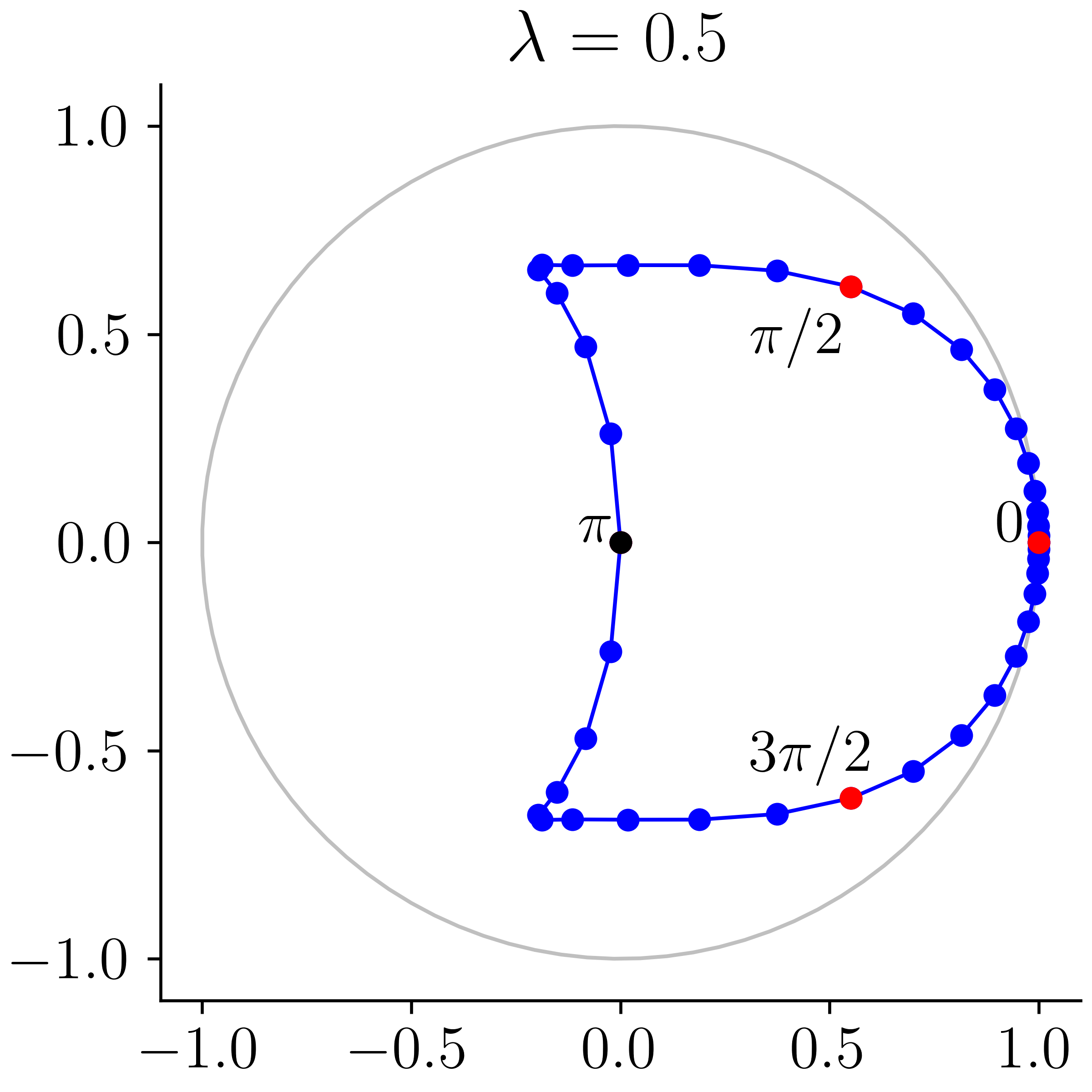}}
\vskip15pt 
\centerline{\includegraphics[trim=0.1cm 0.1cm 0.1cm 0.1cm, clip, width=4.5cm]{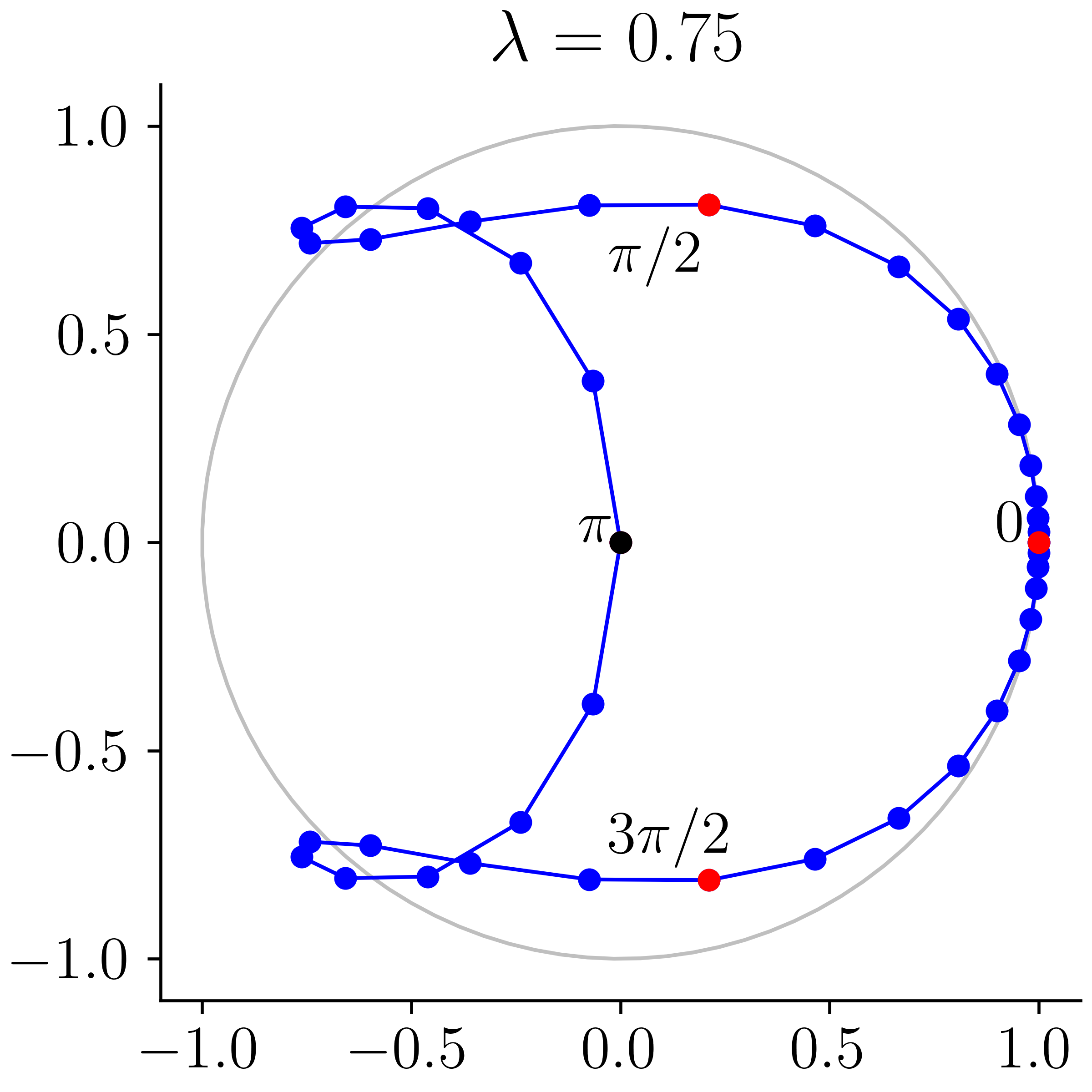}\hspace{1cm}
            \includegraphics[trim=0.1cm 0.1cm 0.1cm 0.1cm, clip, width=4.5cm]{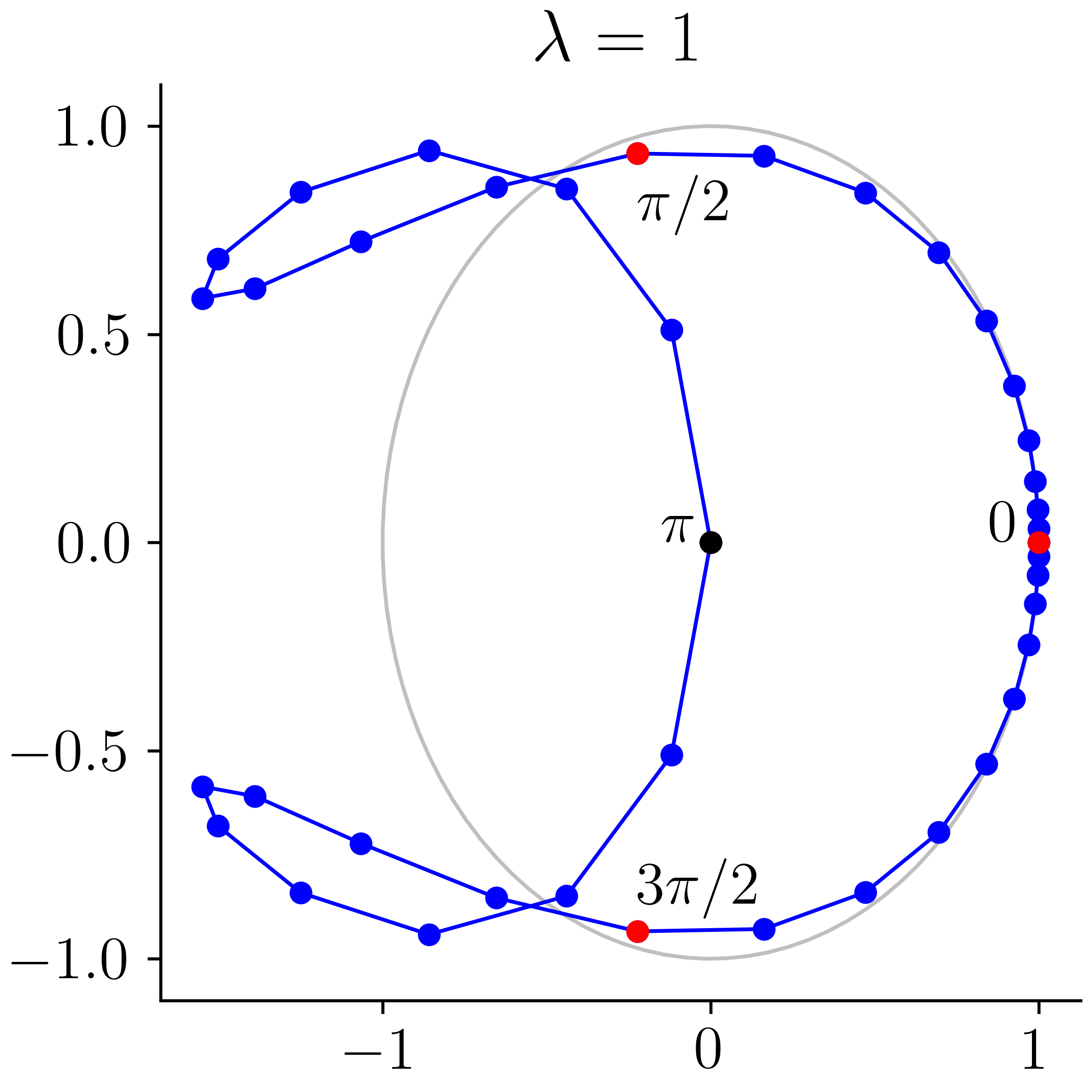}}            
\caption{\sf QL seventh-order scheme: The values $g(\lambda,B_*,\alpha)$ relative to the unit circle in the complex plane for
$\alpha=\alpha_\ell$ given by \eref{21aa} and different $\lambda$.\label{fig5}}
\end{figure}

Note that the value of $B_*$ is independent of $\lambda$ since $\delta(\pi)=0$. It is instructive, however, to study the stability property
of the optimal QL seventh-order scheme for different $\lambda$. To this end, we plot the values of $g(\lambda,B_*,\alpha_\ell)$ for
$\lambda=0.25$, $0.5$, $0.75$, and $1$ in Figure \ref{fig5}, where one can see that the optimal QL seventh-order scheme becomes
unstable if large CFL numbers are used. 

\subsection{Von Neumann Stability Analysis of the Quasi-Linear Ninth-Order Schemes}\label{sec2c}
\begin{figure}[ht!]
\centerline{\includegraphics[trim=0.1cm 0.1cm 0.1cm 0.1cm, clip, width=4.5cm]{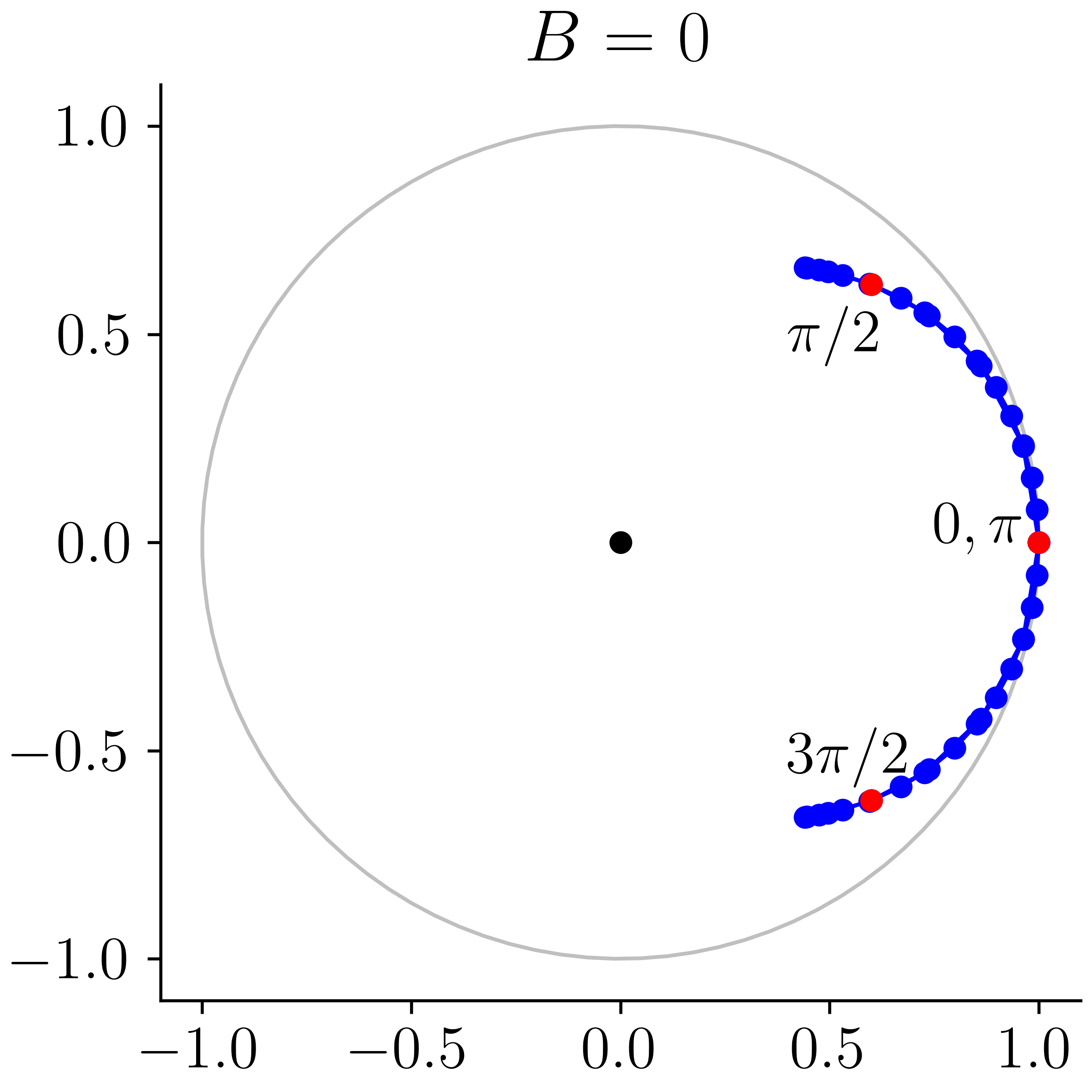}\hspace{1cm}
            \includegraphics[trim=0.1cm 0.1cm 0.1cm 0.1cm, clip, width=4.5cm]{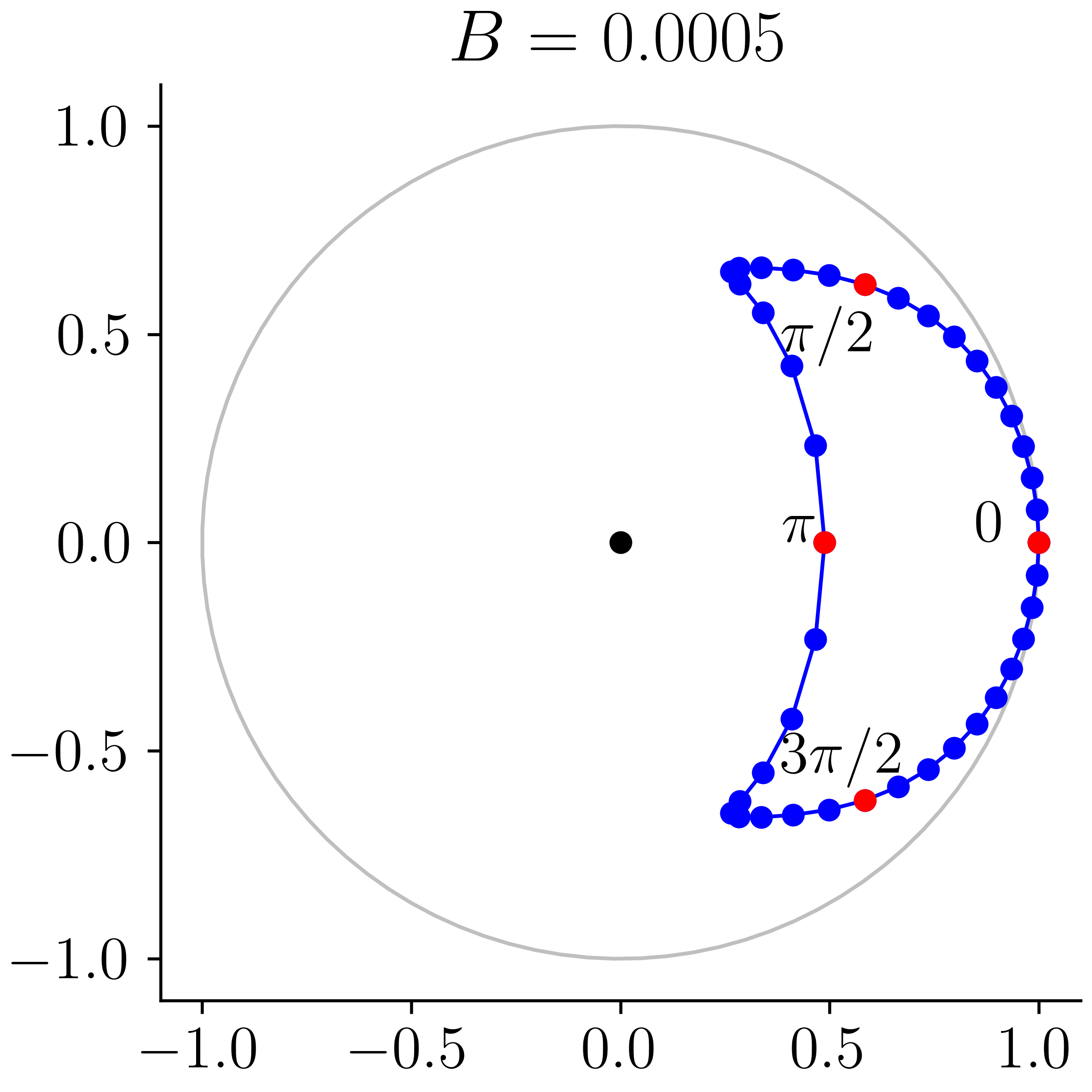}}
\vskip15pt 
\centerline{\includegraphics[trim=0.1cm 0.1cm 0.1cm 0.1cm, clip, width=4.5cm]{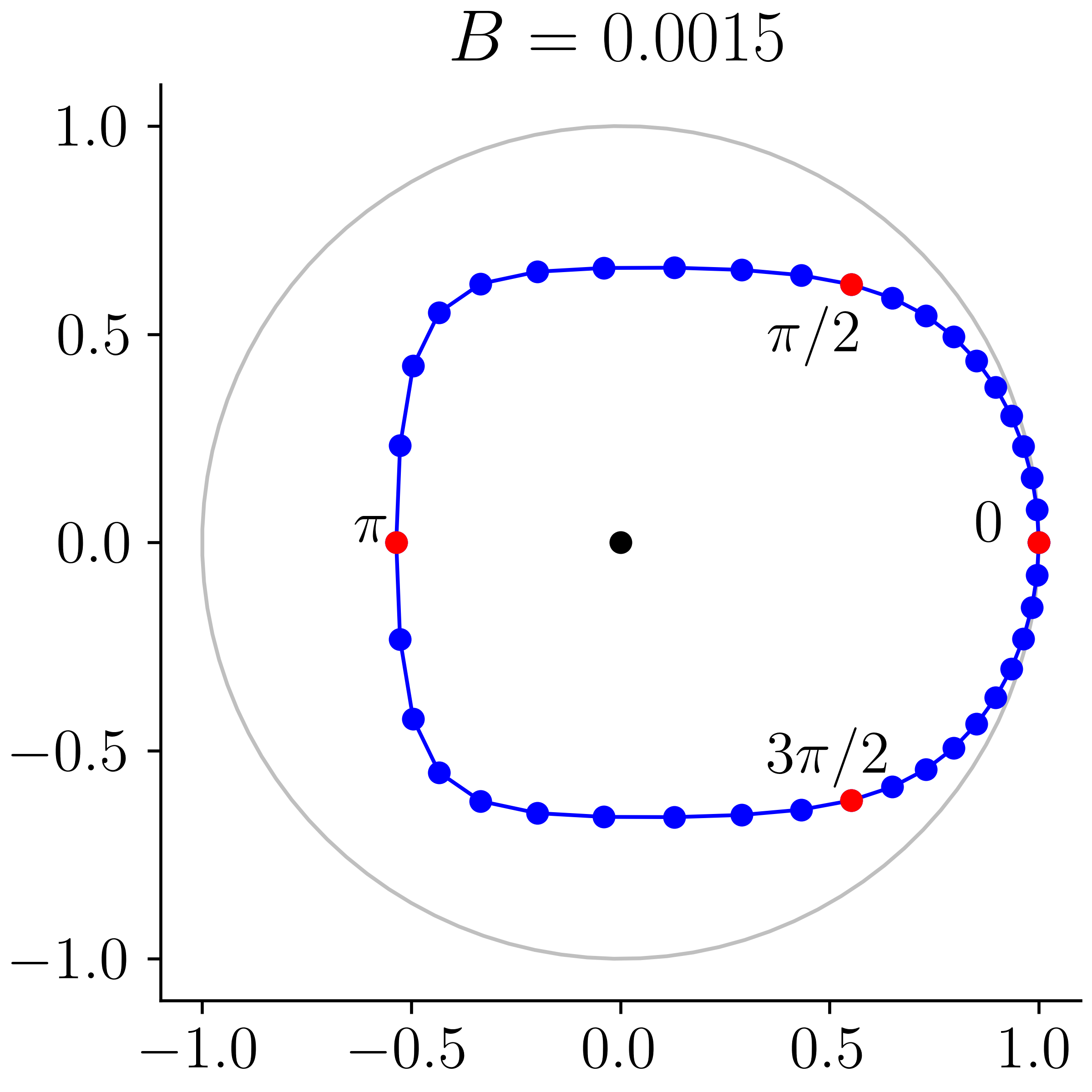}\hspace{1cm}
            \includegraphics[trim=0.1cm 0.1cm 0.1cm 0.1cm, clip, width=4.5cm]{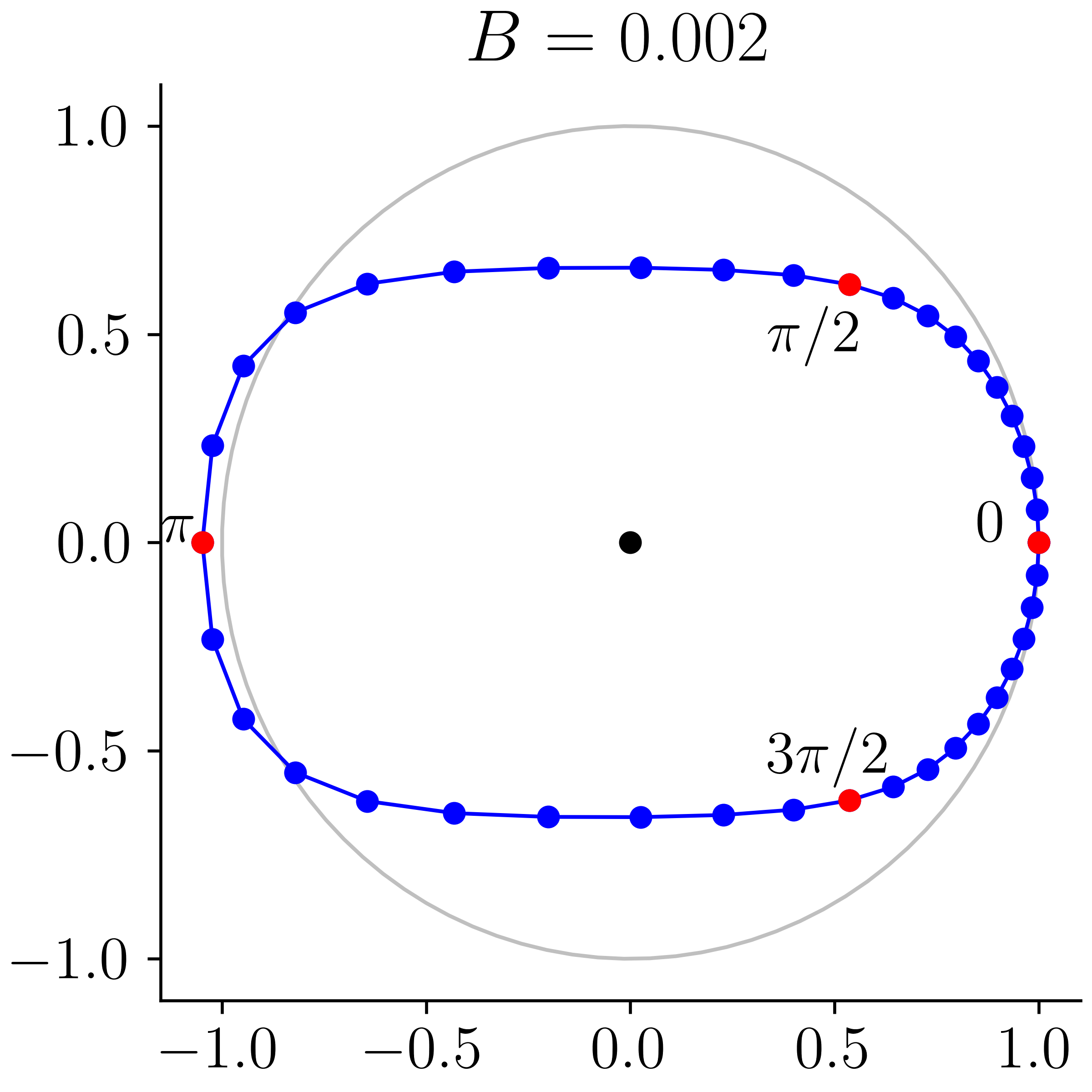}}            
\caption{\sf QL ninth-order scheme: The values $g=g(0.5,B,\alpha)$ relative to the unit circle in the complex plane for $\alpha=\alpha_\ell$
given by \eref{21aa} and different $B$.\label{fig2}}
\end{figure}

We now follow the approach presented in \S \ref{sec2b} and perform the von Neumann stability analysis for the QL ninth-order scheme, which
reads as \eref{2.1bb}--\eref{2.1cc} with
\begin{equation}
\begin{aligned}
&\bm{{\cal L}}_\jph[\bm U]=\frac{1}{2520}\big[2\mF(\mU_{j+5})-23\mF(\mU_{j+4})+127\mF(\mU_{j+3})-473\mF(\mU_{j+2})+1627\mF(\mU_{j+1})\\
&\hspace{2.9cm}+1627\mF(\mU_j)-473\mF(\mU_{j-1})+127\mF(\mU_{j-2})-23\mF(\mU_{j-3})+2\mF(\mU_{j-4})\big],\\
&\bm\omega_\jph^n=
\frac{3}{2048\lambda^n}\big[\mU_{j+5}^n-9\mU_{j+4}^n+36\mU_{j+3}^n-84\mU_{j+2}^n+126\mU_{j+1}^n\\
&\hspace{2.8cm}-126\mU_j^n+84\mU_{j-1}^n-36\mU_{j-2}^n+9\mU_{j-3}^n-\mU_{j-4}^n\big].
\end{aligned}
\label{2.1bb2}
\end{equation}
We apply these schemes to the linear advection equation \eref{Neum1} (once again, in this case, $\lambda\equiv\lambda$), substitute the
Fourier mode \eref{Neum2} into the scheme \eref{2.1bb}--\eref{2.1cc}, \eref{2.1bb2}, and collect like terms, yielding the amplification
factor \eref{Neum6} with 
\begin{equation}
\begin{aligned}
\delta(\alpha)&=-\frac{1}{1260}(e^{5i\alpha}-e^{-5i\alpha})+\frac{5}{504}(e^{4i\alpha}-e^{-4i\alpha})-
\frac{5}{84}(e^{3i\alpha}-e^{-3i\alpha})\\
&\hspace*{0.45cm}+\frac{5}{21}(e^{2i\alpha}-e^{-2i\alpha})-\frac{5}{6}(e^{i\alpha}-e^{-i\alpha}) \\
&=-i\left(\frac{1}{630}\sin(5\alpha)-\frac{5}{252}\sin(4\alpha)+\frac{5}{42}\sin(3\alpha)-\frac{10}{21}\sin(2\alpha)+\frac{5}{3}\sin(\alpha)
\right),\\
\Theta(\alpha)&= e^{5i\alpha}+ e^{-5i\alpha}-10(e^{4i\alpha}+e^{-4i\alpha})+45(e^{3i\alpha}+e^{-3i\alpha})-120(e^{2i\alpha}+e^{-2i\alpha})\\
&\hspace*{0.40cm}+210(e^{i\alpha}+e^{-i\alpha})-252\\
&=2\cos(5\alpha)-20\cos(4\alpha)+90\cos(3\alpha)-240\cos(2\alpha)+420\cos\alpha-252.
\end{aligned}
\label{Neum12}
\end{equation}
In Figure \ref{fig2}, we plot the values of $g$ for $\lambda=0.5$ for different values $\alpha$ defined by \eref{21aa}, and $B=0$, $0.0005$,
$0.0015$, and $0.002$. As one can see, when $B=0.002$, the scheme is unstable, while for $B=0.0005$ and $0.0015$, it is stable. To determine
the optimal value of $B$, we substitute \eref{Neum6} and \eref{Neum12} into \eref{Neum11}, which results in $B_*=1/1024$.

\begin{figure}[ht!]
\centerline{\includegraphics[trim=0.1cm 0.1cm 0.1cm 0.1cm, clip, width=4.5cm]{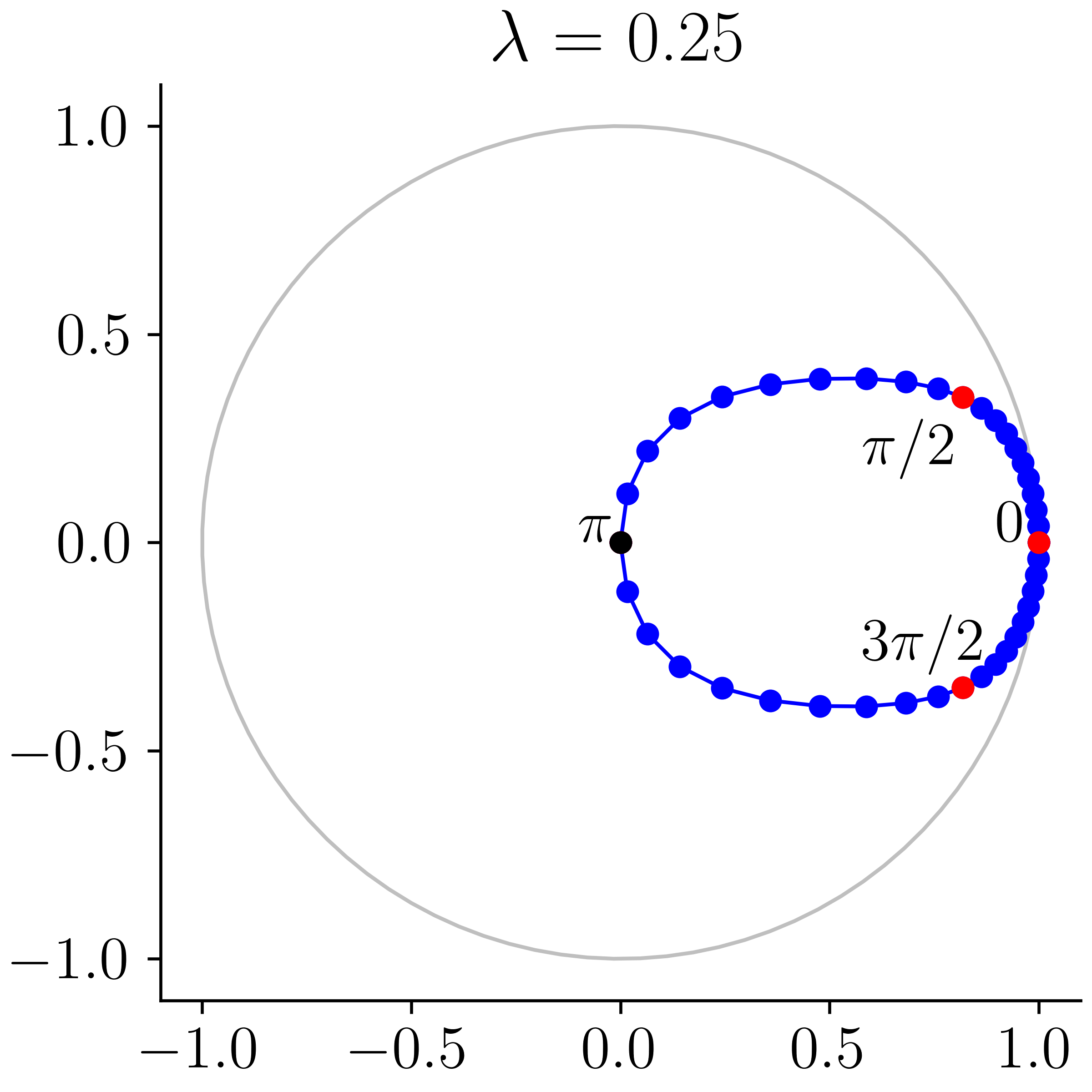}\hspace{1cm}
            \includegraphics[trim=0.1cm 0.1cm 0.1cm 0.1cm, clip, width=4.5cm]{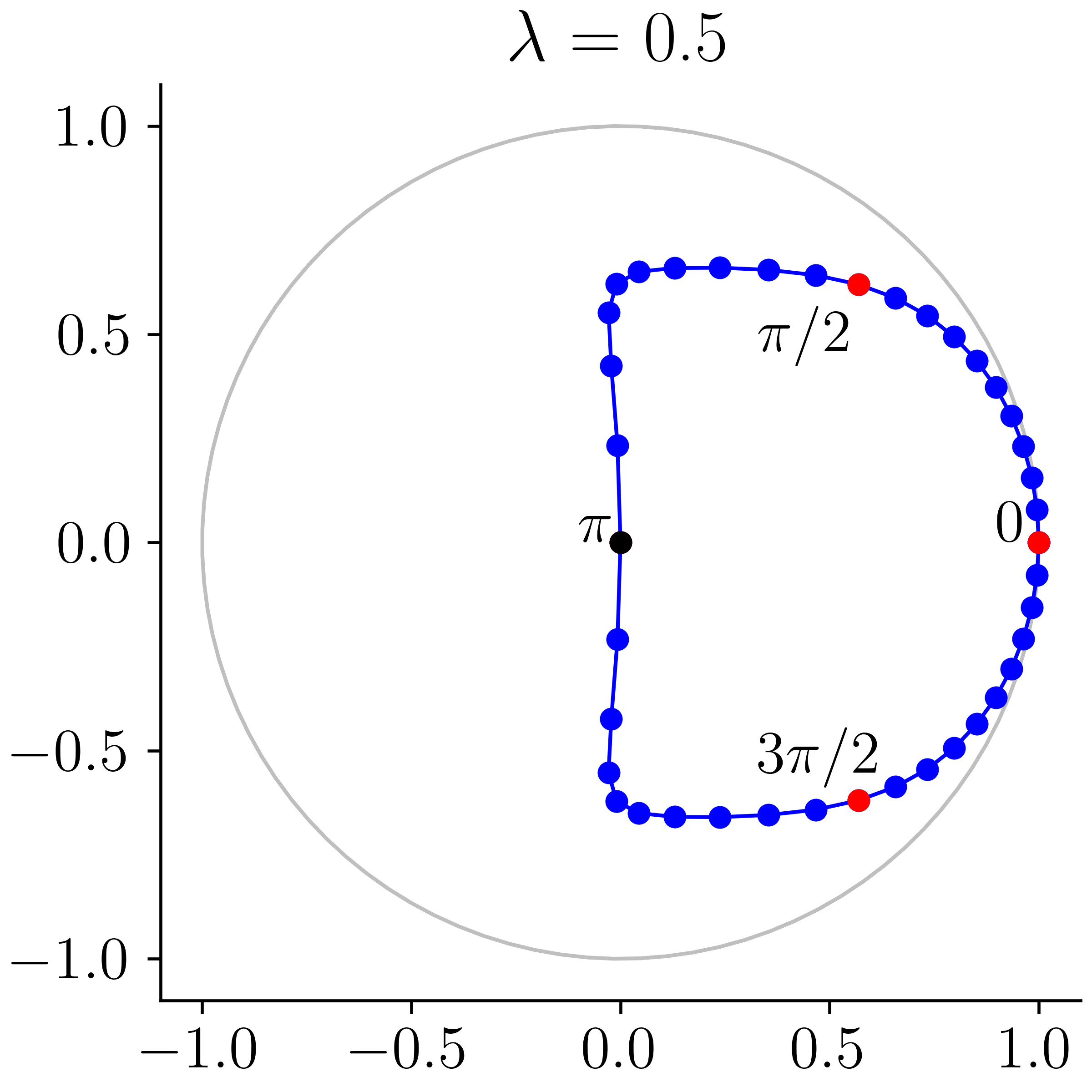}}
\vskip15pt
\centerline{\includegraphics[trim=0.1cm 0.1cm 0.1cm 0.1cm, clip, width=4.5cm]{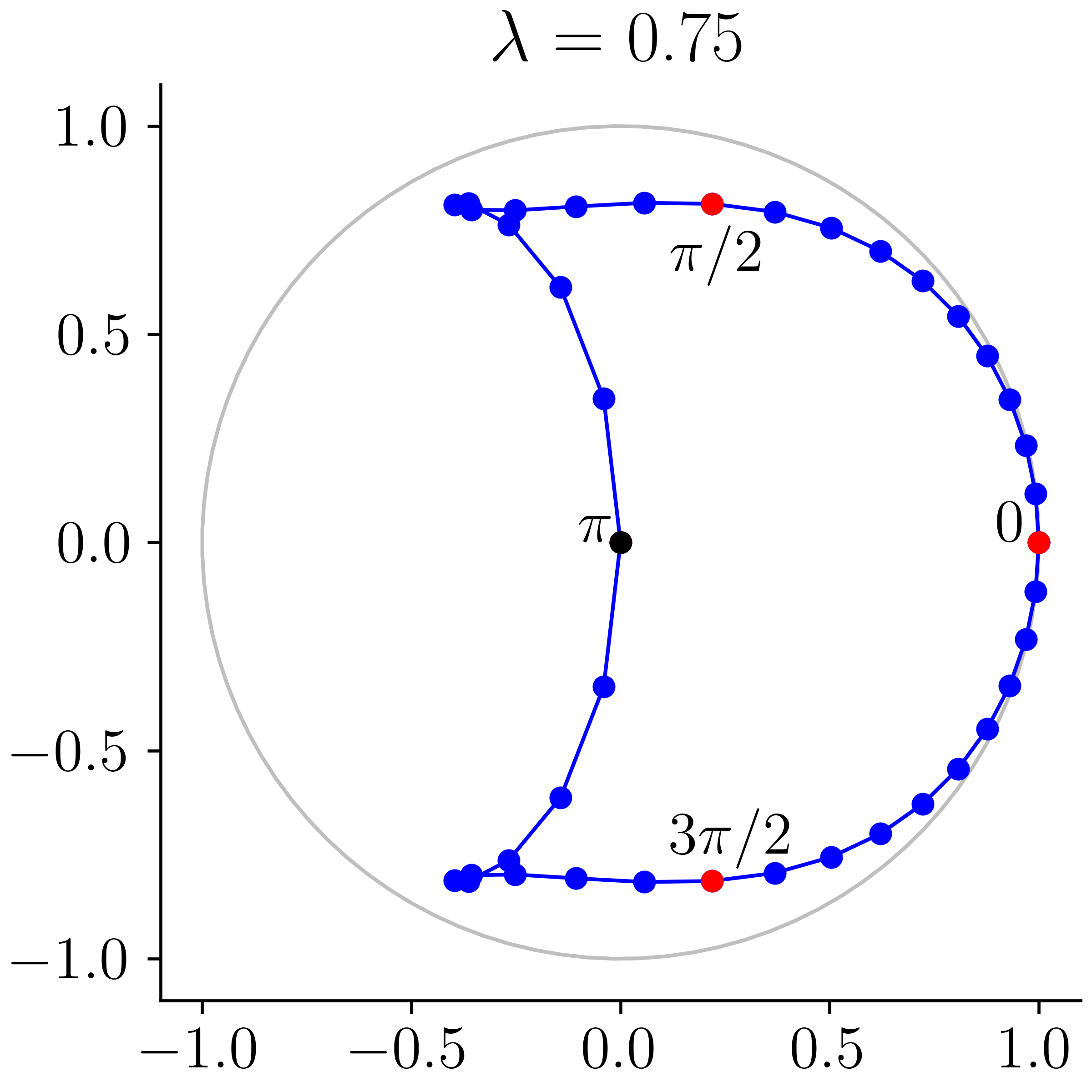}\hspace{1cm}
            \includegraphics[trim=0.1cm 0.1cm 0.1cm 0.1cm, clip, width=4.5cm]{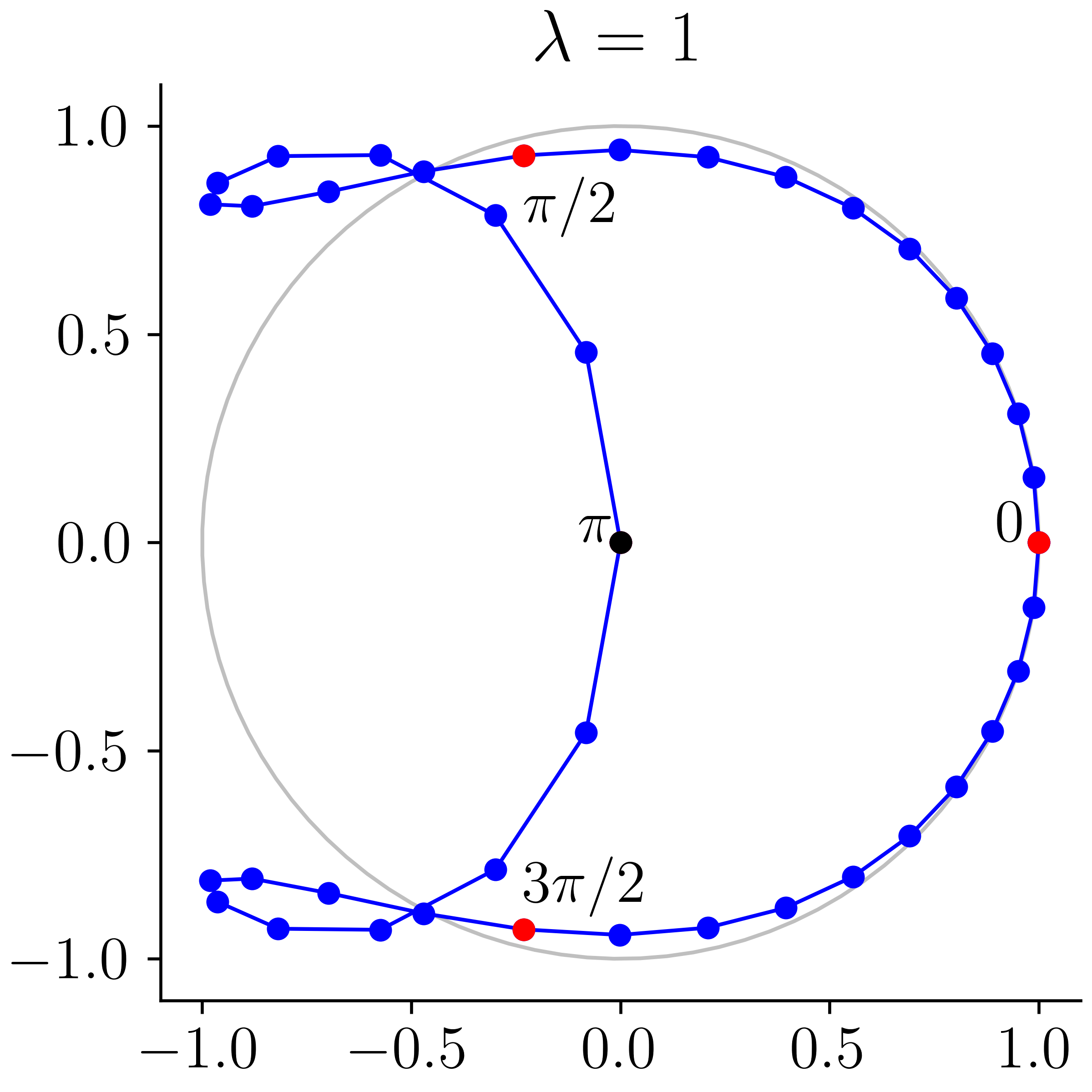}}            
\caption{\sf QL ninth-order scheme: The values $g(\lambda,B_*,\alpha)$ relative to the unit circle in the complex plane for
$\alpha=\alpha_\ell$ given by \eref{21aa} and different $\lambda$.\label{fig6}}
\end{figure}

We also plot the values of $g(\lambda,B_*,\alpha_\ell)$ for different $\lambda$ in Figure \ref{fig6}, where one can see that the optimal QL
ninth-order scheme admits a larger stability region and satisfies the von Neumann condition up to $\lambda=0.75$. 

\subsection{1-D Accuracy Test}
In this section, we demonstrate that the described 1-D QL schemes achieve seventh and ninth orders of accuracy in capturing smooth
solutions. We consider the Euler equations of gas dynamics, which in the 1-D case, read as \eref{1.1} with $\mU=(\rho,\rho u,E)^\top$ and
$\mF=\big(\rho u,\rho u^2+p,u(E+p)\big)^\top$. Here, $\rho$, $u$, $p$, and $E$ are the density, velocity, pressure, and total energy,
respectively, and the system is closed through the following equations of state (EOS) for ideal gases:
\begin{equation*}
p=(\gamma-1)\Big[E-\hf\rho u^2\Big],
\end{equation*}
where the parameter $\gamma$ represents the specific heat ratio. 

We consider the following smooth initial data \cite{KKOKC}:
\begin{equation*}
u(x,0)=\sin\Big(\frac{\pi x}{5}+\frac{\pi}{4}\Big),\quad
\rho(x,0)=\bigg[\frac{\gamma-1}{2\sqrt{\gamma}}\left(u(x,0)+10\right)\bigg]^{\frac{2}{\gamma-1}},\quad p(x,0)=\rho^\gamma(x,0),
\end{equation*}
subject to the periodic boundary conditions in the computational domain $[0,10]$. We compute the numerical solution until the final time
$t=0.1$ on a sequence of uniform meshes with $\dx=1/2$, $1/4$, $1/8$, and $1/16$.

We then compute the $L^1$-errors and estimate the experimental convergence rates for the density $\rho$ using the following Runge formulae,
which are based on the solutions computed on the three consecutive uniform grids with the mesh sizes $\dx$, $2\dx$, and $4\dx$ and denoted
by $(\cdot)^{\dx}$, $(\cdot)^{2\dx}$, and $(\cdot)^{4\dx}$, respectively:
$$
{\rm Error}(\dx)\approx\frac{\delta_{12}^2}{|\delta_{12}-\delta_{24}|},\quad
{\rm Rate}(\dx)\approx\log_2\left(\frac{\delta_{24}}{\delta_{12}}\right).
$$
Here, $\delta_{12}:=\|(\cdot)^{\dx}-(\cdot)^{2\dx}\|_{L^1}$ and $\delta_{24}:=\|(\cdot)^{2\dx}-(\cdot)^{4\dx}\|_{L^1}$. The obtained results
are reported in Table \ref{tab1}, where one can clearly see that the expected orders of accuracy are achieved. 
\begin{table}[ht!]
\centering
\begin{tabular}{ccccccccccccc}
\toprule 
\multirow{2}{1em}{$\dx$}&\multicolumn{2}{c}{QL Seventh-Order Scheme}&\multicolumn{2}{c}{QL Ninth-Order Scheme}\\
\cline{2-5}&Error&Rate&Error&Rate\\
\hline
$1/8$ &2.28e-08&7.05&4.03e-10&8.99\\
$1/16$&1.80e-10&7.02&8.14e-13&8.97\\
\bottomrule 
\end{tabular}
\caption{\sf The $L^1$-errors and experimental convergence rates for the QL seventh- and ninth-order schemes.\label{tab1}}
\end{table}
\begin{rmk}
We stress that in order to achieve the seventh and ninth order of accuracy, we have used smaller time steps with
$\dt\sim(\dx)^{\frac{7}{3}}$ and $\dt\sim(\dx)^3$ to balance the spatial and temporal errors for the seventh- and ninth-order QL schemes,
respectively. 
\end{rmk}

\subsection{1-D Adaptive Schemes}
The 1-D adaptive schemes are similar to the scheme introduced in \cite[\S3]{CFKKO2025}, except that the QL fifth-order scheme from
\cite{KOK25} is now replaced by either the seventh- or ninth-order QL schemes in the smooth regions, which we use with the optimal values
$B=B_*$.

\section{Two-Dimensional Schemes}\label{sec3}
We now extend the 1-D adaptive schemes described in \S\ref{sec2} to the 2-D hyperbolic systems of conservation laws \eref{1.2}. We assume
that the computational domain is covered with uniform cells $I_{j,k}:=\big[x_\jmh,x_\jph\big]\times\big[y_\kmh,y_\kph\big]$ with
$x_\jph-x_\jmh\equiv\dx$ and $y_\kph-y_\kmh\equiv\dy$ centered at $(x_j,y_k)$ with $x_j=\big(x_\jmh+x_\jph\big)/2$ and
$y_k=\big(y_\kmh+y_\kph\big)/2$.

\subsection{Quasi-Linear $r^{\rm th}$-Order Schemes}\label{sec31}
We extend the 1-D QL $r^{\rm th}$-order schemes described in \S\ref{sec21} to the 2-D case, which can be done in a dimension-by-dimension
manner. 

Assuming the point values $\mU^n_{j,k}:\approx\mU(x_j,y_k,t^n)$ are available, the solution is evolved to the next time level $t^{n+1}$
according to 
\begin{equation}
\begin{aligned}
\mU^{\rm I}_{j,k}&=\mU_{j,k}^n-\lambda^n\left(\bm{{\cal F}}_{\jph,k}-\bm{{\cal F}}_{\jmh,k}\right)
-\mu^n\left(\bm{{\cal G}}_{j,\kph}-\bm{{\cal G}}_{j,\kmh}\right),\\
\mU^{\rm II}_{j,k}&=\frac{3}{4}\,\mU_{j,k}^n+\frac{1}{4}\left[\mU^{\rm I}_{j,k}-
\lambda^n\left(\bm{{\cal F}}^{\rm I}_{\jph,k}-\bm{{\cal F}}_{\jmh,k}^{\rm I}\right)\right.
\left.-\mu^n\left(\bm{{\cal G}}^{\rm I}_{j,\kph}-\bm{{\cal G}}^{\rm I}_{j,\kmh}\right)\right],\\[0.5ex]
\mU_{j,k}^{n+1}&=\frac{1}{3}\,\mU_{j,k}^n+\frac{2}{3}\left[\mU^{\rm II}_{j,k}-
\lambda^n\left(\bm{{\cal F}}^{\rm II}_{\jph,k}-\bm{{\cal F}}^{\rm II}_{\jmh,k}\right)\right.
\left.-\mu^n\left(\bm{{\cal G}}^{\rm II}_{j,\kph}-\bm{{\cal G}}^{\rm II}_{j,\kmh}\right)\right].
\end{aligned}
\label{3.1}
\end{equation}
where $\lambda^n:=\dt^n/\dx$, $\mu^n:=\dt^n/\dy$, and the numerical fluxes are defined by
\begin{equation}
\begin{aligned}
&\bm{{\cal F}}_{\jph,k}=\bm{{\cal L}}_{\jph,k}^x[\bm U^n],&&
\bm{{\cal F}}^{\rm I}_{\jph,k}=\bm{{\cal L}}_{\jph,k}^x\big[\bm U^{\rm I}\big],&&
\bm{{\cal F}}^{\rm II}_{\jph,k}=\bm{{\cal L}}_{\jph,k}^x\big[\bm U^{\rm II}\big]-\bm\omega^{x,n}_{\jph,k},\\
&\bm{{\cal G}}_{j,\kph}=\bm{{\cal L}}_{j,\kph}^y[\bm U^n],&&
\bm{{\cal G}}_{j,\kph}^{\rm I}=\bm{{\cal L}}_{j,\kph}^y\big[\bm U^{\rm I}\big],&&
\bm{{\cal G}}_{j,\kph}^{\rm II}=\bm{{\cal L}}_{j,\kph}^y\big[\bm U^{\rm II}\big]-\bm\omega^{y,n}_{j,\kph}.
\end{aligned}
\label{3.2}
\end{equation}
Here, 
\begin{equation*}
\begin{aligned}
&\bm{{\cal L}}^x_{\jph,k}[\bm U]:=\sum_{i=1-\ell}^\ell b_i\mF(\mU_{j+i,k}),
&&\bm\omega_{\jph,k}^{x,n}:=\frac{3B_*}{2\lambda^n}\sum_{i=1-\ell}^\ell(-1)^{i+\ell}\binom{2\ell-1}{i+\ell-1}\mF(\mU_{j+i,k}),\\
&\bm{{\cal L}}^y_{j,\kph}[\bm U]:=\sum_{i=1-\ell}^\ell b_i\mG(\mU_{j,k+i}),
&&\bm\omega_{j,\kph}^{y,n}:=\frac{3B_*}{2\mu^n}\sum_{i=1-\ell}^\ell(-1)^{i+\ell}\binom{2\ell-1}{i+\ell-1}\mG(\mU_{j,k+i}), 
\end{aligned}
\end{equation*}
where $b_i$ are given by \eref{2.1dd} and the optimal value of $B_*$ is the same as the 1-D case, namely, $B_*=1/{2^{r+1}}$.

\subsection{2-D Accuracy Test}
In this section,  we demonstrate that the developed 2-D QL schemes achieve the seventh and ninth orders of accuracy in capturing smooth
solutions. As in the 1-D case, we focus on the Euler equations of gas dynamics, whose 2-D version reads as \eref{1.2} with
\begin{equation}
\mU=(\rho,\rho u,\rho v,E)^\top,\quad\mF=\big(\rho u,\rho u^2+p,\rho uv,u(E+p)\big)^\top,\quad
\mG=\big(\rho v,\rho uv,\rho v^2+p,v(E+p)\big)^\top.
\label{3.3}
\end{equation}
Here, $v$ is the $y$-velocity and the other variables are the same as in the 1-D case. The system \eref{1.2}, \eref{3.3} is closed through
the following EOS for ideal gases:
\begin{equation}
p=(\gamma-1)\Big[E-\frac{\rho}{2}(u^2+v^2)\Big].
\label{3.4}
\end{equation}

We consider the following smooth initial data \cite{BD2013,Shu1998}:
\begin{equation*}
\begin{aligned}
&\rho(x,y,0)=\bigg(1-\frac{(\gamma-1)\kappa^2}{2\gamma}\bigg)^{\frac{1}{\gamma-1}},\quad p(x,y,0)=\rho^\gamma(x,y,0),\\
&u(x,y,0)=1-\kappa y,\quad v(x,y,0)=1+\kappa x,\quad\kappa=\frac{5}{2\pi}e^{\frac{1-x^2-y^2}{2}},
\end{aligned}
\end{equation*}
subject to the periodic boundary conditions in the computational domain $[-10,10]\times[-10,10]$. The exact solution of this initial value
problem is given by $\mU(x,y,t)=\mU(x-t,y-t,0)$.

We compute the numerical solution until the final time $t=0.1$ on a sequence of uniform meshes with $\dx=\dy=1/10$, $1/20$, and $1/40$. The
time steps are selected to be $\dt\sim(\dx)^{\frac{7}{3}}$ and $\dt\sim(\dx)^3$ for the seventh- and ninth-order schemes, respectively. We
then measure the $L^1$-errors and the corresponding experimental convergence rates for the density $\rho$. The obtained results are
presented in Table \ref{tab52}, where one can see that the expected orders of accuracy have been achieved.
\begin{table}[ht!]
\centering
\begin{tabular}{ccccccccccccc}
\toprule 
\multirow{2}{5em}{$\dx=\dy$}&\multicolumn{2}{c}{QL Seventh-Order Scheme}&\multicolumn{2}{c}{QL Ninth-Order Scheme}\\
\cline{2-5}&Error&Rate&Error&Rate \\
\hline
$1/10$&2.01e-05&--- &7.09e-06&---\\
$1/20$&1.57e-07&7.00&1.43e-08&8.96\\
$1/40$&1.27e-09&6.95&2.88e-11&8.95\\
\bottomrule 
\end{tabular}
\caption{\sf The $L^1$-errors and experimental convergence rates for the QL seventh- and ninth-order schemes.\label{tab52}}
\end{table}

\subsection{2-D Adaptive Schemes}
As in the 1-D case, the 2-D adaptive schemes are similar to the one introduced in \cite[\S4.3]{CFKKO2025}, except that the QL fifth-order
schemes from \cite{KOK25} are replaced by either the seventh- or ninth-order QL schemes in the smooth regions. These schemes read as
\eref{3.1}--\eref{3.2} with 
\begin{equation*}
\begin{aligned}
&\bm{{\cal L}}_{\jph,k}^x[\bm U]=\frac{1}{840}\big[-3\mF(\mU_{j+4,k})+29\mF(\mU_{j+3,k})-139\mF(\mU_{j+2,k})+533\mF(\mU_{j+1,k})\\
&\hspace{3.05cm}+533\mF(\mU_{j,k})-139\mF(\mU_{j-1,k})+29\mF(\mU_{j-2,k})-3\mF(\mU_{j-3,k})\big],\\
&\bm\omega_{\jph,k}^{x,n}=\frac{3}{512\lambda^n}\big[-\mU_{j+4,k}^n+7\mU_{j+3,k}^n-21\mU_{j+2,k}^n+35\mU_{j+1,k}^n-
35\mU_{j,k}^n+21\mU_{j-1,k}^n-7\mU_{j-2,k}^n+\mU_{j-3,k}^n\big],\\[0.5ex]
&\bm{{\cal L}}_{j,\kph}^y[\bm U]=\frac{1}{840}\big[-3\mG(\mU_{j,k+4})+29\mG(\mU_{j,k+3})-139\mG(\mU_{j,k+2})+533\mG(\mU_{j,k+1})\\
&\hspace{3.05cm}+533\mG(\mU_{j,k})-139\mG(\mU_{j,k-1})+29\mG(\mU_{j,k-2})-3\mG(\mU_{j,k-3})\big],\\
&\bm\omega_{j,\kph}^{y,n}=\frac{3}{512\mu^n}\big[-\mU_{j,k+4}^n+7\mU_{j,k+3}^n-21\mU_{j,k+2}^n+35\mU_{j,k+1}^n-
35\mU_{j,k}^n+21\mU_{j,k-1}^n-7\mU_{j,k-2}^n+\mU_{j,k-3}^n\big],
\end{aligned}
\end{equation*}
and 
\begin{equation*}
\begin{aligned}
&\bm{{\cal L}}_{\jph,k}^x[\bm U]=\frac{1}{2520}\big[2\mF(\mU_{j+5,k})-23\mF(\mU_{j+4,k})+127\mF(\mU_{j+3,k})-473\mF(\mU_{j+2,k})+
1627\mF(\mU_{j+1,k})\\
&\hspace{3.1cm}+1627\mF(\mU_{j,k})-473\mF(\mU_{j-1,k})+127\mF(\mU_{j-2,k})-23\mF(\mU_{j-3,k})+2\mF(\mU_{j-4,k})\big],\\
&\bm\omega_{\jph,k}^{x,n}=\frac{3}{2048\lambda^n}\big[\mU_{j+5,k}^n-9\mU_{j+4,k}^n+36\mU_{j+3,k}^n-84\mU_{j+2,k}^n+126\mU_{j+1,k}^n\\
&\hspace{3.0cm}-126\mU_{j,k}^n+84\mU_{j-1,k}^n-36\mU_{j-2,k}^n+9\mU_{j-3,k}^n-\mU_{j-4,k}^n\big],\\[0.5ex]
&\bm{{\cal L}}_{j,\kph}^y[\bm U]=\frac{1}{2520}\big[2\mG(\mU_{j,k+5})-23\mG(\mU_{j,k+4})+127\mG(\mU_{j,k+3})-473\mG(\mU_{j,k+2})+
1627\mG(\mU_{j,k+1})\\
&\hspace{3.1cm}+1627\mG(\mU_{j,k})-473\mG(\mU_{j,k-1})+127\mG(\mU_{j,k-2})-23\mG(\mU_{j,k-3})+2\mG(\mU_{j,k-4})\big],\\
&\bm\omega_{j,\kph}^{y,n}=\frac{3}{2048\mu^n}\big[\mU_{j,k+5}^n-9\mU_{j,k+4}^n+36\mU_{j,k+3}^n-84\mU_{j,k+2}^n+126\mU_{j,k+1}^n\\
&\hspace{3.0cm}-126\mU_{j,k}^n+84\mU_{j,k-1}^n-36\mU_{j,k-2}^n+9\mU_{j,k-3}^n-\mU_{j,k-4}^n\big],
\end{aligned}
\end{equation*}
for the seventh- and ninth-order schemes, respectively.

\section{Numerical Examples}\label{sec4}
In this section, we apply the developed adaptive schemes to both the 1-D and 2-D Euler equations of gas dynamics and compare their
performance with that of the adaptive schemes from \cite{CFKKO2025}. We take $\gamma=1.4$ in Examples 1--5 and $\gamma=5/3$ in Example 6. In
all of the examples, we use the CFL number $0.4$.

For the sake of brevity, the adaptive schemes which are based on the QL fifth-, seventh-, and ninth-order schemes will be referred to as the
5-Order, 7-Order, and 9-Order schemes, respectively.

In the examples below, we specify the values of the adaption constants $\texttt{C}_1$ and $\texttt{C}_2$. For their specific definitions, we
refer the reader to \cite{CFKKO2025}, where details on the adaptive algorithms can be found. 

All of the numerical experiments have been conducted on a workstation equipped with an Intel(R) Core(TM) i9-13900H CPU at 2.6 GHz and 32 GB
of RAM. The simulations were conducted in FORTRAN using GCC version 14.2.0 compiler suite. The reported CPU times were averaged over thirty
independent runs to ensure reproducibility and minimize variability due to system processes.

\subsection{1-D Examples}
In the 1-D examples, the adaptive solutions will be plotted along with the reference solutions, which will be computed using the
second-order LDCU scheme from \cite{CKX_24} and implemented with the Minmod2 limiter on much finer meshes. 

\paragraph{Example 1---Titarev-Toro Problem.} In the first example taken from \cite{Toro2005} (see also \cite{Shu88,Toro2005a}), we consider
the shock-entropy wave interaction problem with the following initial conditions:
\begin{equation*}
(\rho,u,p)\Big|_{(x,0)}=\begin{cases}(1.51695,0.523346,1.805),&x<-4.5,\\(1+0.1\sin(20x),0,1),&x>-4.5,\end{cases}
\end{equation*}
prescribed in the computational domain $[-5,5]$ subject to the free boundary conditions. In this initial boundary value problem, a
forward-facing shock wave of Mach number $1.1$ interacts with high-frequency density perturbations: as the shock wave moves, the
perturbations spread ahead.

\begin{figure}[ht!]
\centerline{\includegraphics[trim=0.9cm 0.4cm 0.9cm 0.6cm, clip, width=7.0cm]{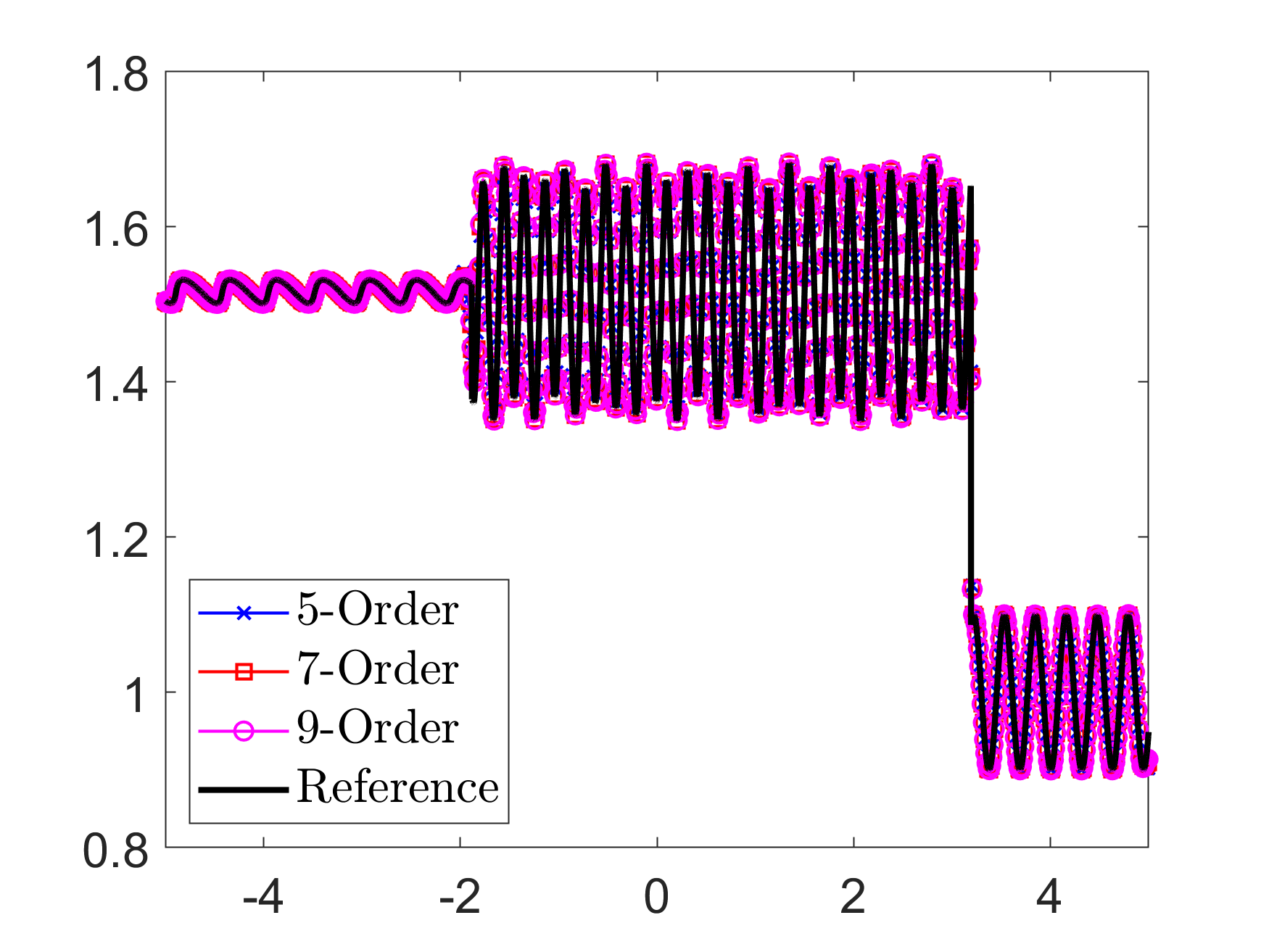}\hspace{1cm}
            \includegraphics[trim=0.9cm 0.4cm 0.9cm 0.6cm, clip, width=7.0cm]{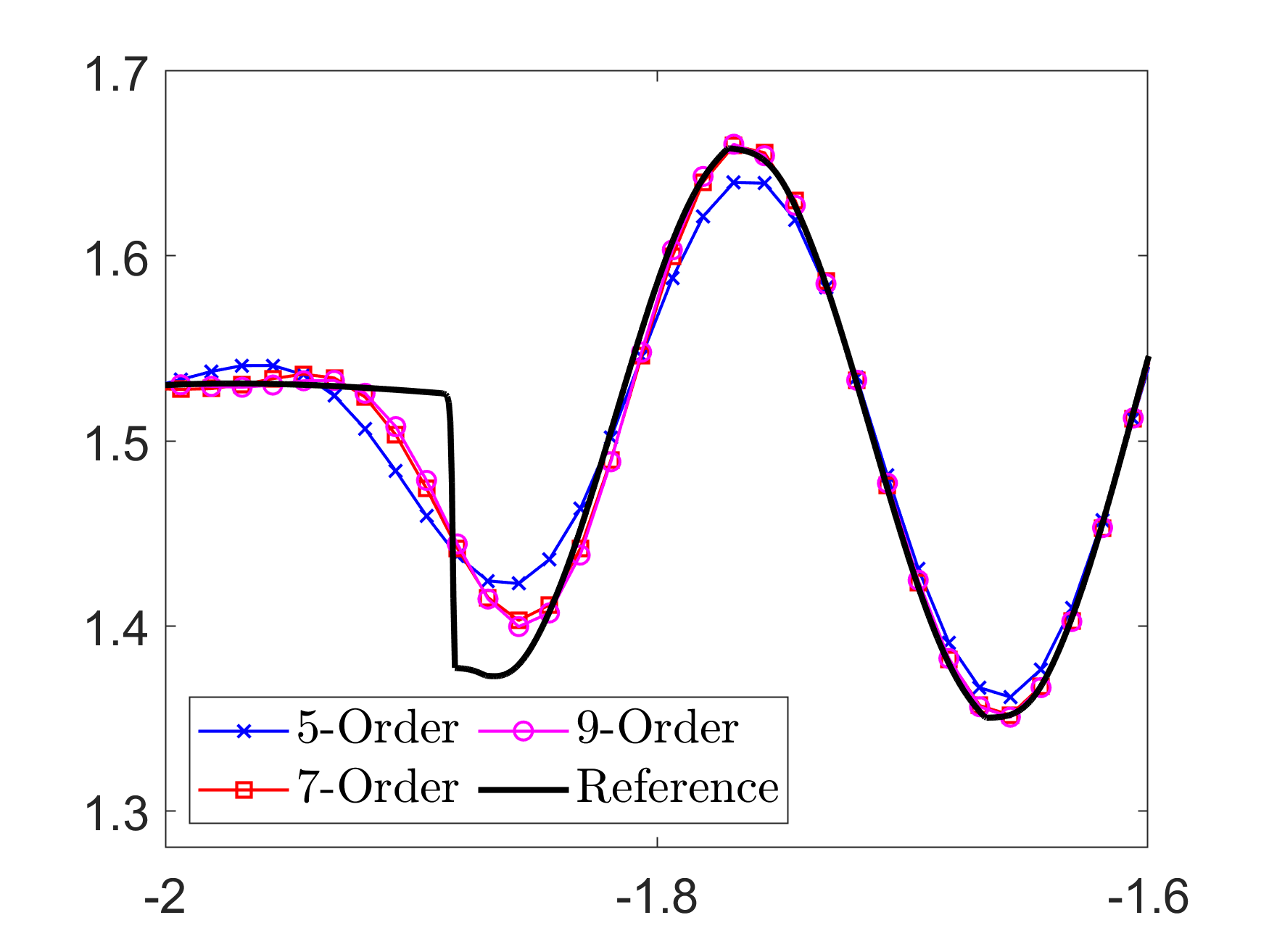}}
\caption{\sf Example 1: Density $\rho$ computed by the 5-Order, 7-Order, and 9-Order schemes (left) and zoom at $x\in[-2,-1.6]$ (right).
\label{fig3a}}
\end{figure}

We compute the numerical solutions by the 5-Order, 7-Order, and 9-Order schemes (with the adaption constants $\texttt{C}_1=0.02$ and
$\texttt{C}_2=0.3$) until the final time $t=5$ on a uniform mesh with $\dx=1/60$. The obtained numerical results are presented in Figure
\ref{fig3a} along with the reference solution computed on a much finer mesh with $\dx=1/1600$, where one can see that the 7-Order and
9-Order schemes are very close but both are more accurate than the 5-Order one.

\paragraph*{Example 2---Blast Wave Problem.} In the second 1-D example, we consider a strong shocks interaction problem from
\cite{Woodward88} with the initial data,
\begin{equation*}
(\rho,u,p)\Big|_{(x,0)}=\begin{cases}(1,0,1000),&x<0.1,\\(1,0,0.01),&0.1\le x\le0.9,\\(1,0,100),&x>0.9,\end{cases}
\end{equation*}
prescribed in the computational domain $[0,1]$ subject to the solid wall boundary conditions.

\begin{figure}[ht!]
\centerline{\includegraphics[trim=1.3cm 0.4cm 0.8cm 0.6cm, clip, width=7.0cm]{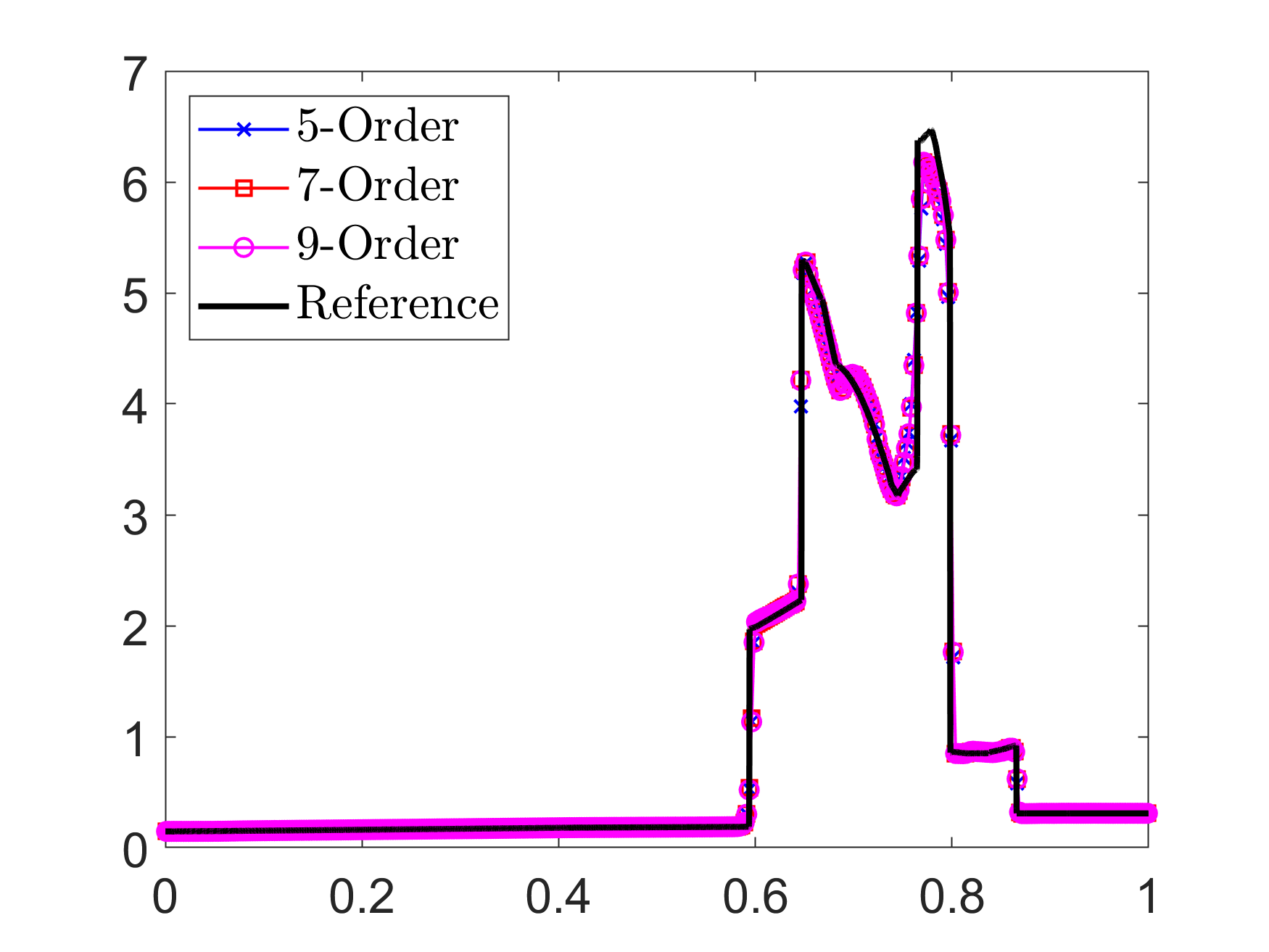}\hspace*{1.0cm}
            \includegraphics[trim=1.3cm 0.4cm 0.8cm 0.6cm, clip, width=7.0cm]{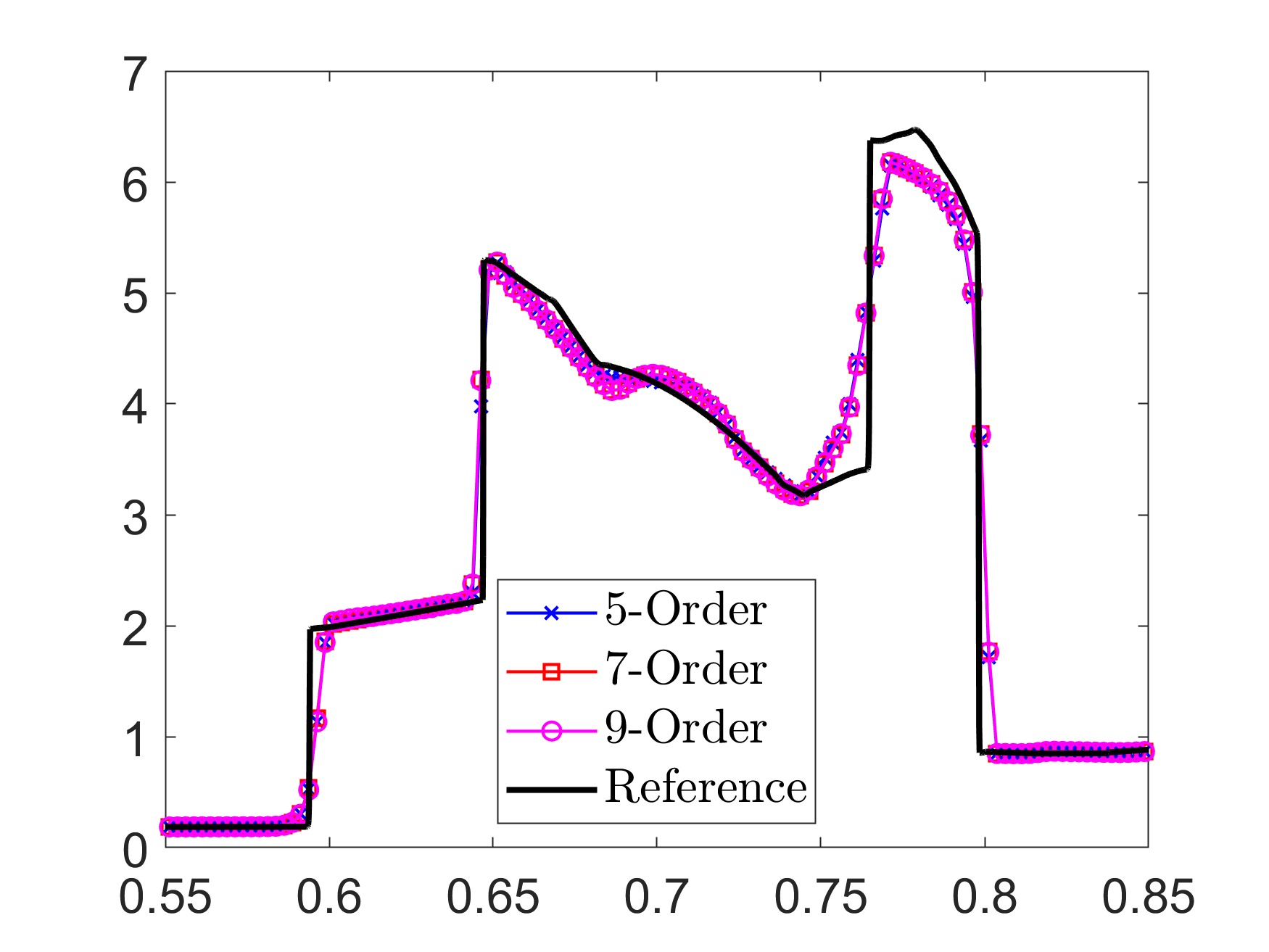}}
\caption{\sf Example 2: Density $\rho$ computed by the 5-Order, 7-Order, and 9-Order schemes (left) and zoom at $x\in[0.55,0.85]$.
\label{fig5a}}
\end{figure}

We compute the numerical solution until the final time $t=5$ by the 5-Order, 7-Order, and 9-Order schemes (with the adaption constants
$\texttt{C}_1=0.002$ and $\texttt{C}_2=0.3$) on a uniform mesh with $\dx=1/400$ and present the obtained numerical results in Figure
\ref{fig5a} together with the reference solution computed on a much finer mesh with $\dx=1/4000$. It can be observed that all of the studied
adaptive schemes are capable of achieving a superb resolution of the contact wave located at about $x=0.6$. The obtained results are quite
similar, but the main goal of this example is to show that the new adaptive schemes do not produce any spurious oscillations while capturing
the numerical solutions of this challenging benchmark. 

\subsection{2-D Examples}
In the 2-D examples, we will not only show the obtained solutions, but will also plot the areas, which were identified as ``rough'' at the
last time step. The contact discontinuities areas will be marked in red, while the rest of the ``rough'' areas will be indicated in black. 

\paragraph*{Example 3---2-D Riemann Problem (Configuration 3).} In this example taken from \cite{Kurganov02}, the initial conditions,
\begin{equation*}
(\rho,u,v,p)\Big|_{(x,y,0)}=\begin{cases}(1.5,0,0,1.5),&x>1,~y>1,\\(0.5323,1.206,0,0.3),&x<1,~y>1,\\(0.138,1.206,1.206,0.029),&x<1,~y<1,\\
(0.5323,0,1.206,0.3),&x>1,~y<1,\end{cases}
\end{equation*}
are prescribed in the computational domain $[0,1.2]\times[0,1.2]$ subject to the free boundary conditions.

\begin{figure}[ht!]
\centerline{\includegraphics[trim=0cm 0cm 0cm 0cm, clip, width=\linewidth]{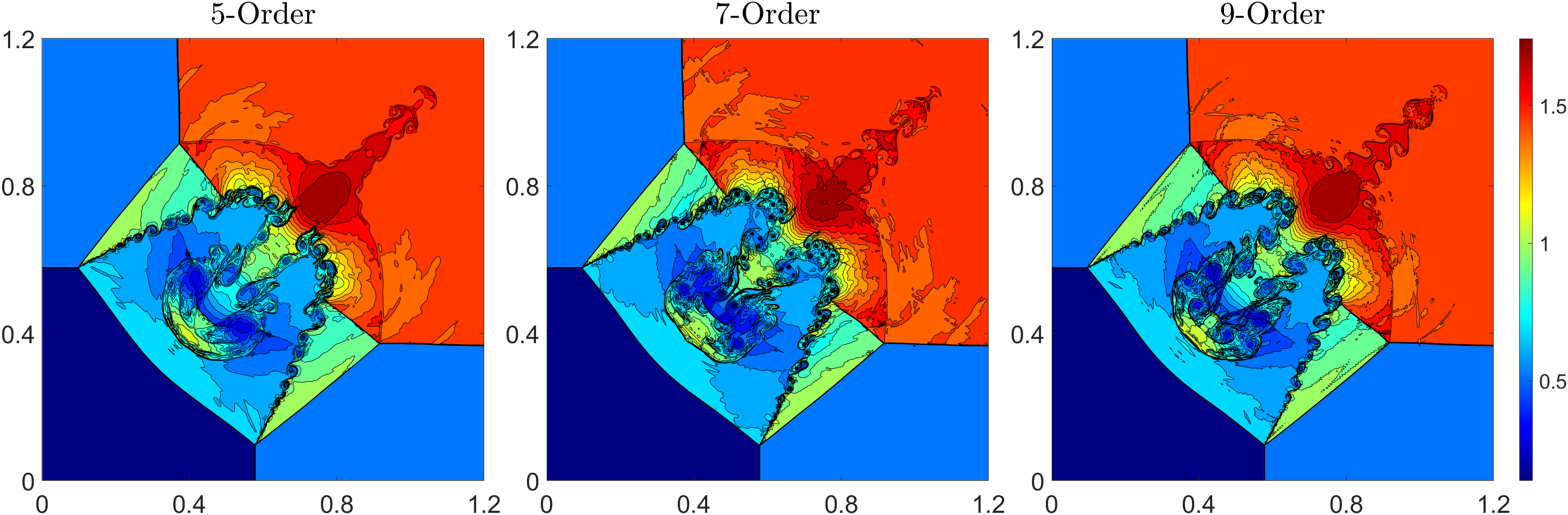}}
\vskip 12pt 
\centerline{\hspace{-0.7cm}\includegraphics[trim=0cm 0cm 0cm 0cm, clip, width=16.8cm]{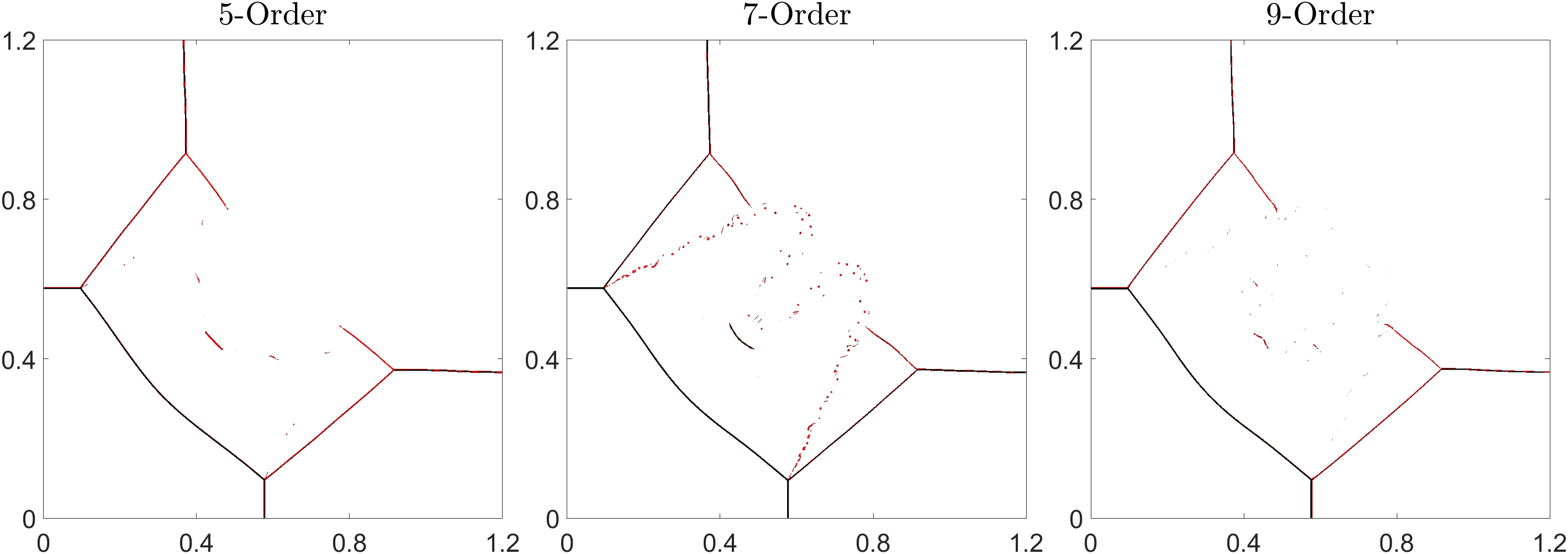}}
\caption{\sf Example 3: Top row: Density $\rho$ computed by the 5-Order (left), 7-Order (middle), and 9-Order (right) schemes. Bottom row:
The ``rough'' areas detected by the 5-Order (left), 7-Order (middle), and 9-Order (right) schemes at the final time step.\label{fig6a}}
\end{figure}

We compute the numerical solution until the final time $t=1$ by the 5-Order (with the adaption constants $\texttt{C}_1=0.08$ and
$\texttt{C}_2=0.2$), 7-Order (with the adaption constants $\texttt{C}_1=0.08$ and $\texttt{C}_2=0.15$), and 9-Order (with the adaption
constants $\texttt{C}_1=0.1$ and $\texttt{C}_2=0.2$) schemes on a uniform mesh with $\dx=\dy=3/2500$. The obtained results are plotted in
Figure \ref{fig6a} (top row), where one can see that the 7-Order and 9-Order schemes slightly outperform the 5-Order one in capturing a
sideband instability of the jet in the zones of strong along-jet velocity shear and the instability along the jets neck. The detected
``rough'' areas for the 5-Order, 7-Order, and 9-Order schemes are plotted in  Figure \ref{fig6a} (bottom row), where one can see that the
overcompressive and Minmod2 limiters are implemented only in a small part of the computational domain.

\paragraph*{Example 4---2-D Riemann Problem (Configuration 6).} In this example also taken from \cite{Kurganov02}, the initial conditions,
\begin{equation*}
(\rho,u,v,p)\Big|_{(x,y,0)}=\begin{cases}(1,0.75,-0.5,1),&x>0.5,~y>0.5,\\(2,0.75,0.5,1),&x<0.5,~y>0.5,\\(1,-0.75,0.5,1),&x<0.5,~y<0.5,\\
(3,-0.75,-0.5,1),&x>0.5,~y<0.5,\end{cases}
\end{equation*}
are prescribed in the computational domain $[0,1]\times[0,1]$ subject to the free boundary conditions.

\begin{figure}[ht!]
\centerline{\includegraphics[trim=0cm 0cm 0cm 0cm, clip, width=\linewidth]{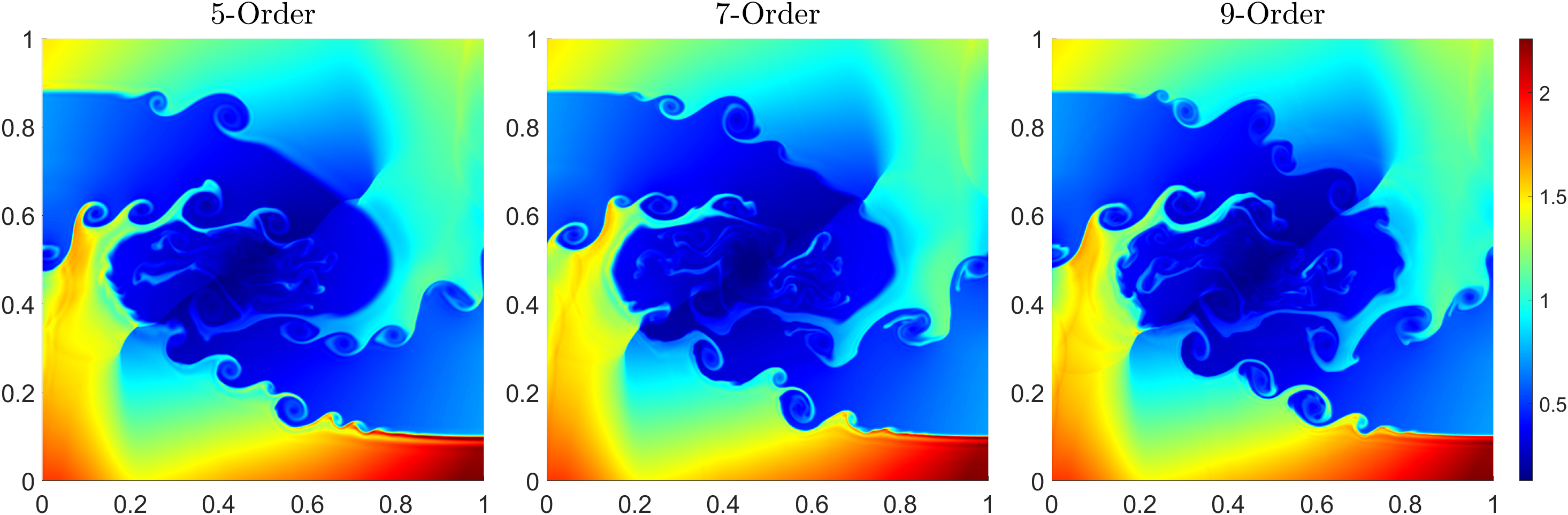}}
\vskip 12pt 
\centerline{\hspace{-0.7cm}\includegraphics[trim=0cm 0cm 0cm 0cm, clip, width=16.8cm]{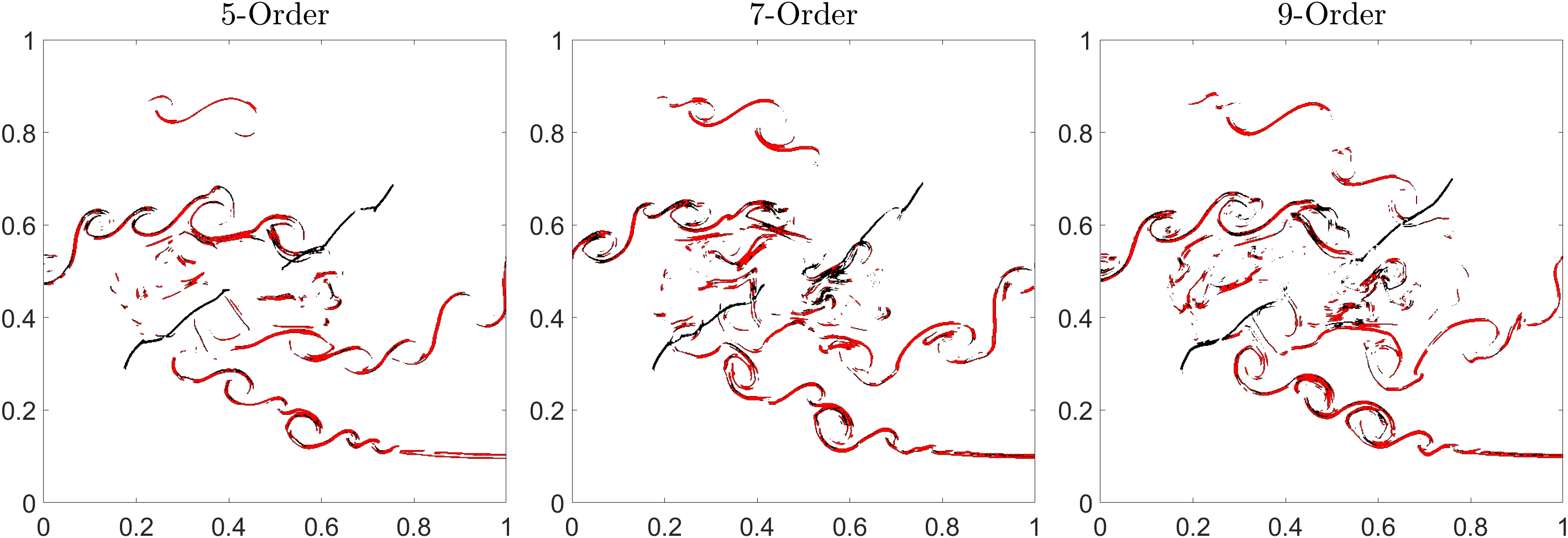}}
\caption{\sf Example 4: Top row: Density $\rho$ computed by the 5-Order (left), 7-Order (middle), and 9-Order (right) schemes. Bottom row:
The ``rough'' areas detected by the 5-Order (left), 7-Order (middle), and 9-Order (right) schemes at the final time step.\label{fig16a}}
\end{figure}

We compute the numerical solution until the final time $t=1$ by the 5-Order, 7-Order, and 9-Order (with the adaption constants
$\texttt{C}_1=0.06$ and $\texttt{C}_2=0.004$) schemes on a uniform mesh with $\dx=\dy=1/600$. The obtained results are plotted in Figure
\ref{fig16a} (top row), where one can see that the 7-Order and 9-Order schemes outperform the 5-Order one in capturing much more complicated
vortex structures and the 9-Order scheme achieves the highest resolution of small solution structures. At the same time, the ``rough'' areas
identified by the 5-Order, 7-Order, and 9-Order schemes are quite similar; see Figure \ref{fig16a} (bottom row).

\paragraph{Example 5---Implosion Problem.} In this example, we consider the implosion problem taken from \cite{Liska03}. The initial
conditions,
\begin{equation*}
(\rho,u,v,p)\Big|_{(x,y,0)}=\begin{cases}(0.125,0,0,0.14),&|x|+|y|<0.15,\\(1,0,0,1),&\mbox{otherwise},\end{cases}
\end{equation*}
are prescribed in the computational domain $[0,0.5]\times[0,0.5]$ with the solid wall boundary conditions. In this setting, a jet forms near
the origin and propagates along the diagonal $y=x$ direction. The position of the jet can serve as an indicator of the amount of numerical
dissipation present in the studied numerical scheme as schemes containing larger numerical dissipation may not resolve the jet at all or the
jet propagation may slow down by excessive numerical diffusion.

\begin{figure}[ht!]
\centerline{\includegraphics[trim=0cm 0cm 0cm 0cm, clip, width=\linewidth]{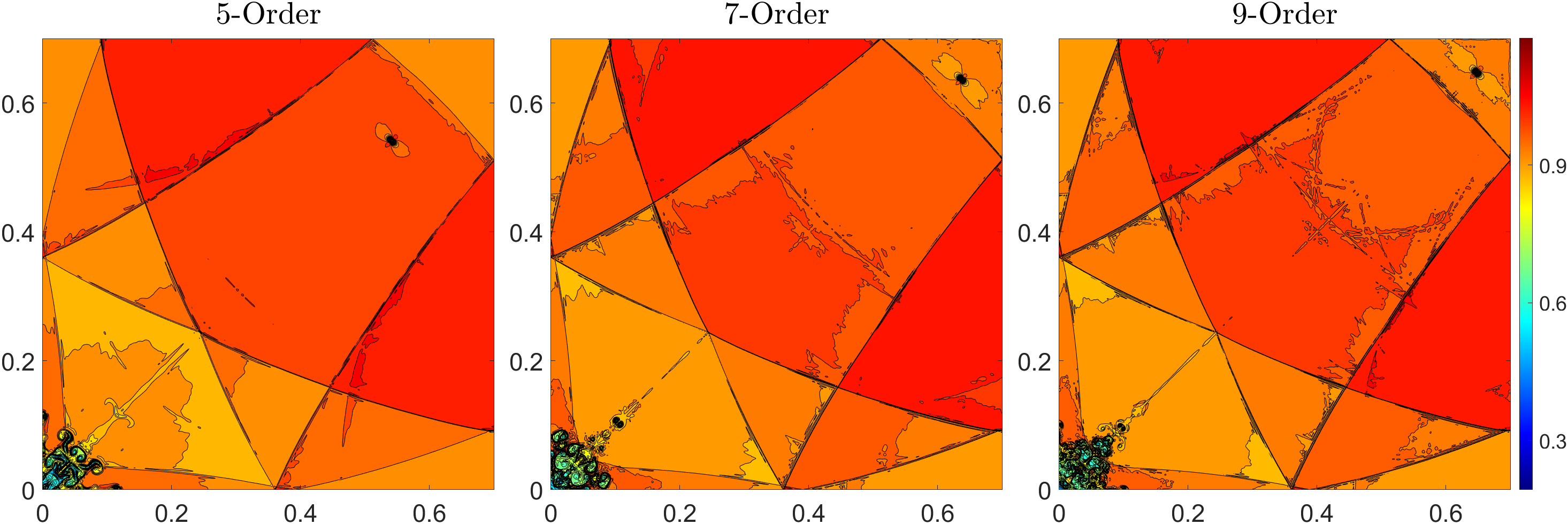}}
\vskip 12pt 
\centerline{\hspace{-0.7cm}\includegraphics[trim=0cm 0cm 0cm 0cm, clip, width=16.8cm]{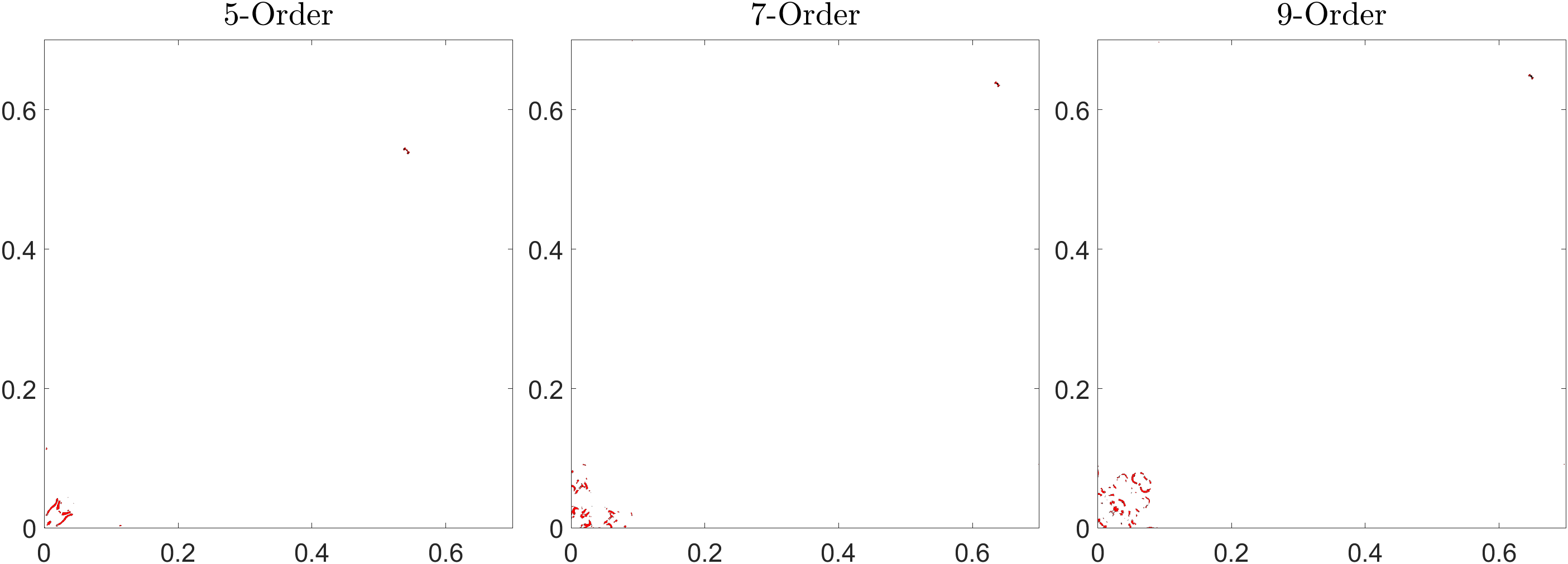}}
\caption{\sf Example 5: Top row: Density $\rho$ computed by the 5-Order (left), 7-Order (middle), and 9-Order (right) schemes. Bottom row:
The ``rough'' areas detected by the 5-Order (left), 7-Order (middle), and 9-Order (right) schemes at the final time step.\label{fig8a}}
\end{figure}

We compute the numerical solution until the final time $t=2.5$ by the 5-Order, 7-Order, and 9-Order (with the adaption constants
$\texttt{C}_1=0.05$ and $\texttt{C}_2=0.1$) schemes on a uniform mesh with $\dx=\dy=7/10000$. We present the obtained numerical results in
Figure \ref{fig8a} (top row), where one can see that the jet propagates further in the diagonal direction in the numerical result computed
by the 7-Order and 9-Order schemes, which demonstrates that the use of higher-order schemes in the smooth areas helps to reduce numerical
dissipation. At the same time, the ``rough'' areas detected by the 5-Order, 7-Order, and 9-Order schemes, which are plotted in Figure
\ref{fig8a} (bottom row), show that nonlinear limiters are implemented only in a very small part of the computational domain.

Next, we measure the CPU times consumed by the 5-Order, 7-Order, and 9-Order schemes. The obtained results show that the CPU time consumed
by the 7-Order and 9-Order schemes are about $7.0\%$ and $11.7\%$ larger than the CPU time consumed by the 5-Order scheme, respectively.
This demonstrates that the proposed adaptive algorithms based on new QL seventh-order and ninth-order schemes achieve higher resolution
compared with their counterpart based on the QL fifth-order scheme, while consuming almost the same CPU times. 

\paragraph{Example 6---RT Instability.} In the last example taken from \cite{Shi03}, we investigate the RT instability, which is a physical
phenomenon occurring when a layer of heavier fluid is placed on top of a layer of lighter fluid.

We consider the 2-D Euler equations of gas dynamics \eref{1.2}, \eref{3.3}--\eref{3.4} with the gravitational source terms acting in the
positive direction of the $y$-axis:
\begin{equation*}
\begin{aligned}
&\rho_t+(\rho u)_x+(\rho v)_y=0,\\
&(\rho u)_t+(\rho u^2 +p)_x+(\rho uv)_y=0,\\
&(\rho v)_t+(\rho uv)_x+(\rho v^2+p)_y=\rho,\\
&E_t+\left[u(E+p)\right]_x+\left[v(E+p)\right]_y=\rho v,
\end{aligned}
\end{equation*}
and take the following initial conditions:
\begin{equation*}
(\rho,u,v,p)\Big|_{(x,y,0)}=
\begin{cases}(2,0,-0.025c\cos(8\pi x),2y+1),&y<0.5,\\(1,0,-0.025c\cos(8\pi x),y+1.5),&\mbox{otherwise},\end{cases}
\end{equation*}
where $c:=\sqrt{\gamma p/\rho}$ is the speed of sound. The solid wall boundary conditions are imposed at $x=0$ and $x=0.25$, and the
following Dirichlet boundary conditions are prescribed at the top and bottom boundaries:
$$
(\rho,u,v,p)(x,1,t)=(1,0,0,2.5),\quad(\rho,u,v,p)(x,0,t)=(2,0,0,1).
$$

\begin{figure}[ht!]
\centerline{\includegraphics[trim=0cm 0cm 0cm 0cm, clip, width=\linewidth]{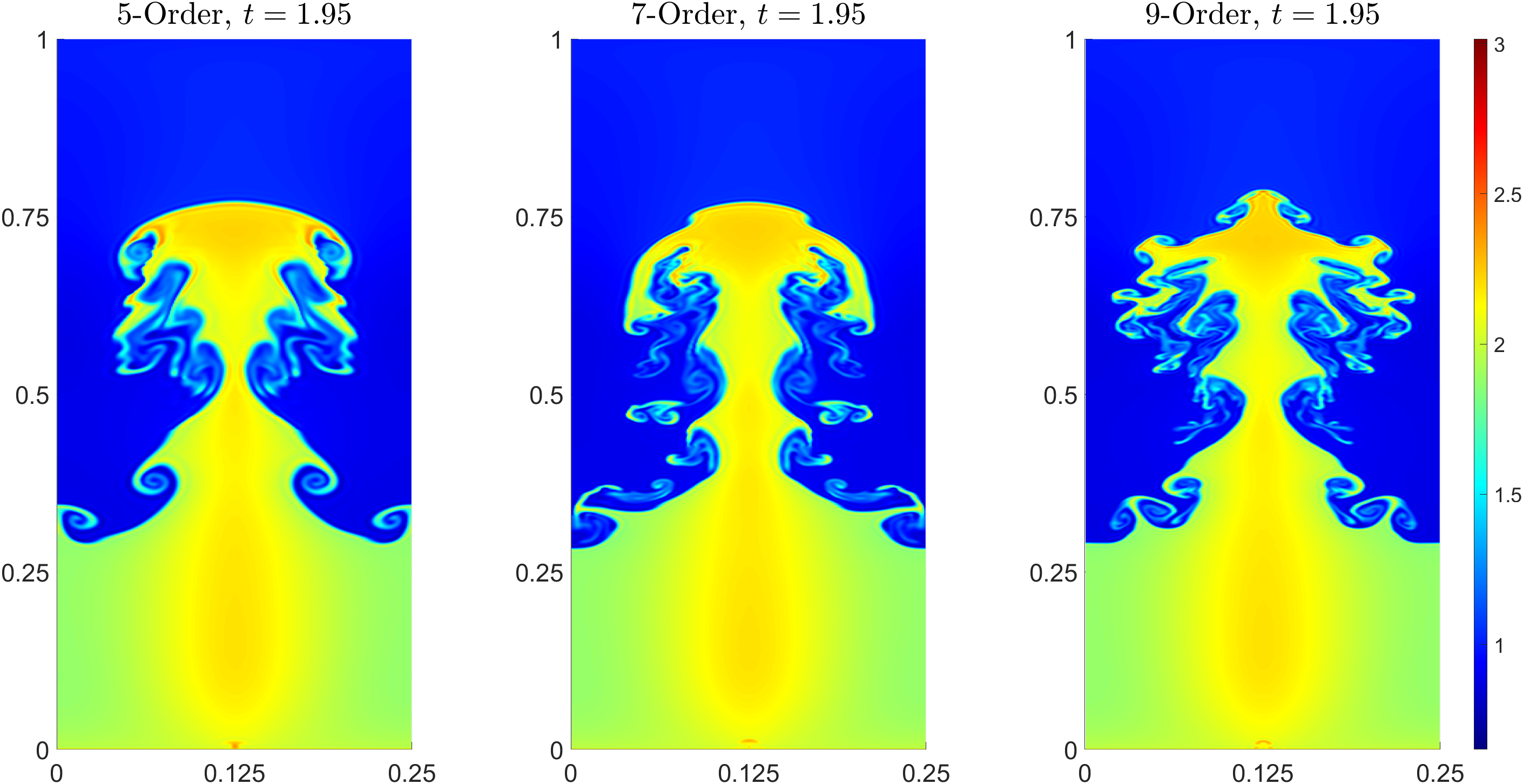}}
\vskip 12pt 
\centerline{\hspace{-0.7cm}\includegraphics[trim=0cm 0cm 0cm 0cm, clip, width=16.8cm]{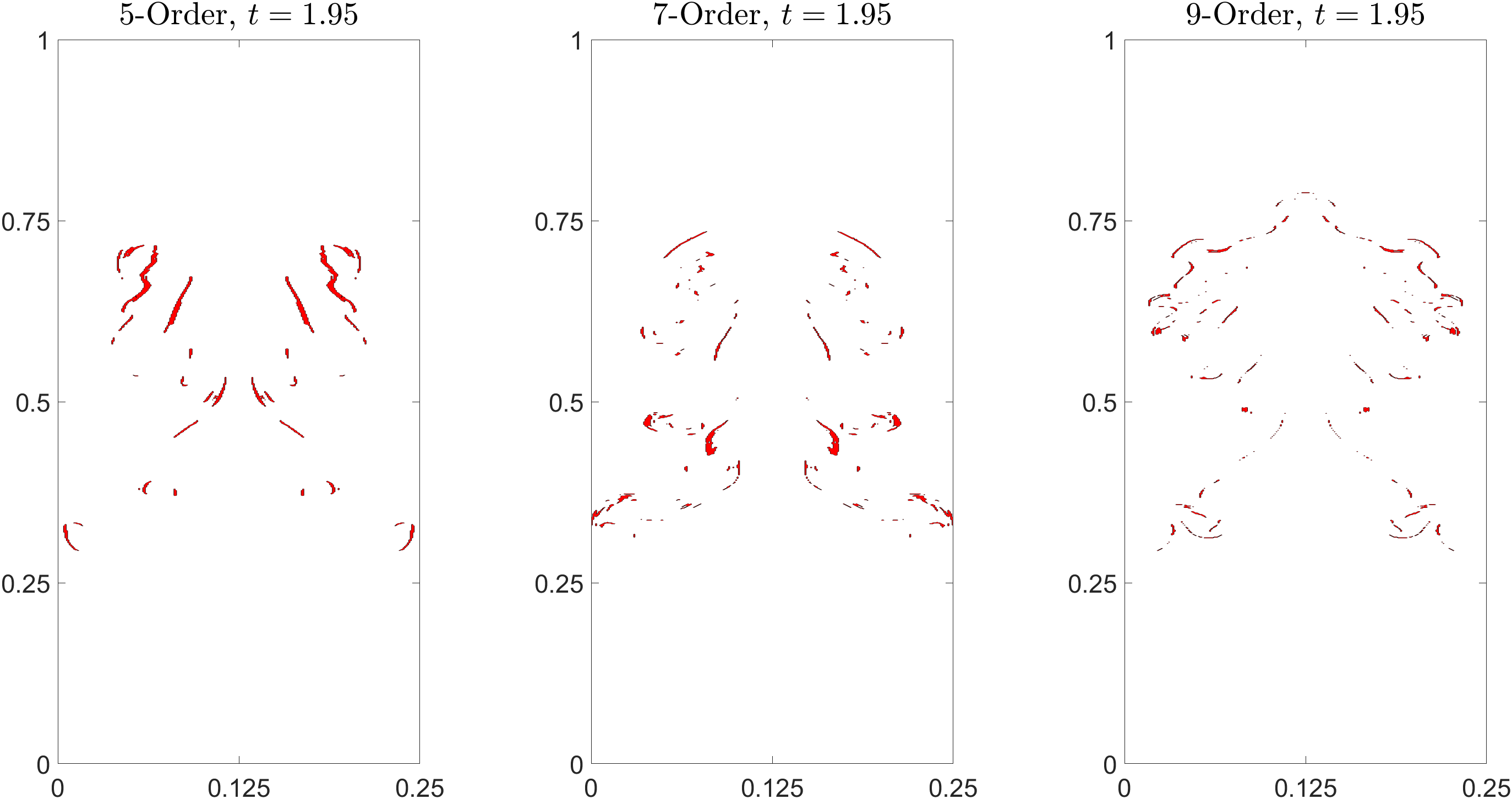}}
\caption{\sf Example 6: Top row: Density $\rho$ computed by the 5-Order (left), 7-Order (middle), and 9-Order (right) schemes. Bottom row:
The ``rough'' areas detected by the 5-Order (left), 7-Order (middle), and 9-Order (right) schemes at $t=1.95$.\label{fig10a}}
\end{figure}

\begin{figure}[ht!] 
\centerline{\includegraphics[trim=0cm 0cm 0cm 0cm, clip, width=\linewidth]{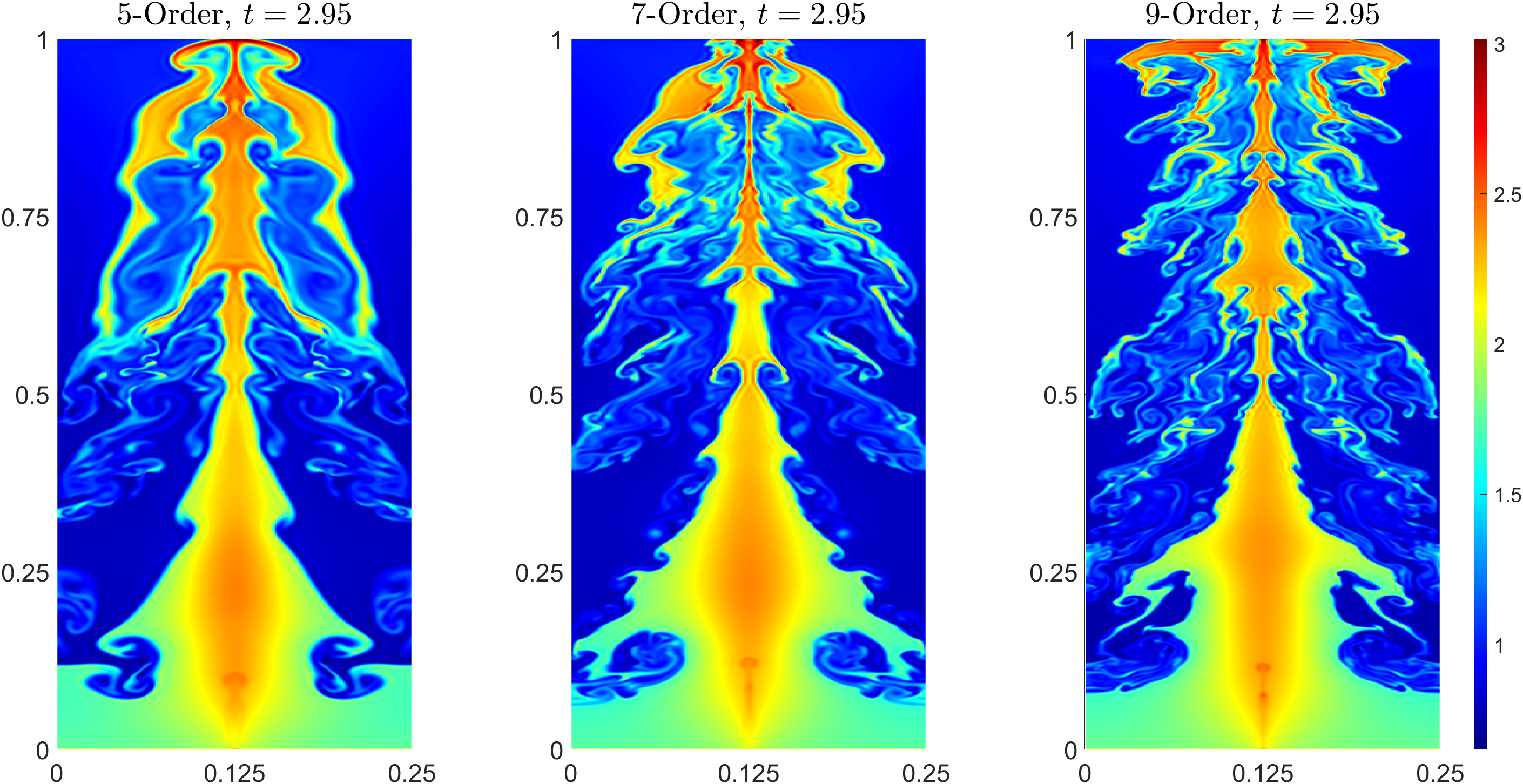}}
\vskip 12pt 
\centerline{\hspace{-0.7cm}\includegraphics[trim=0cm 0cm 0cm 0cm, clip, width=16.8cm]{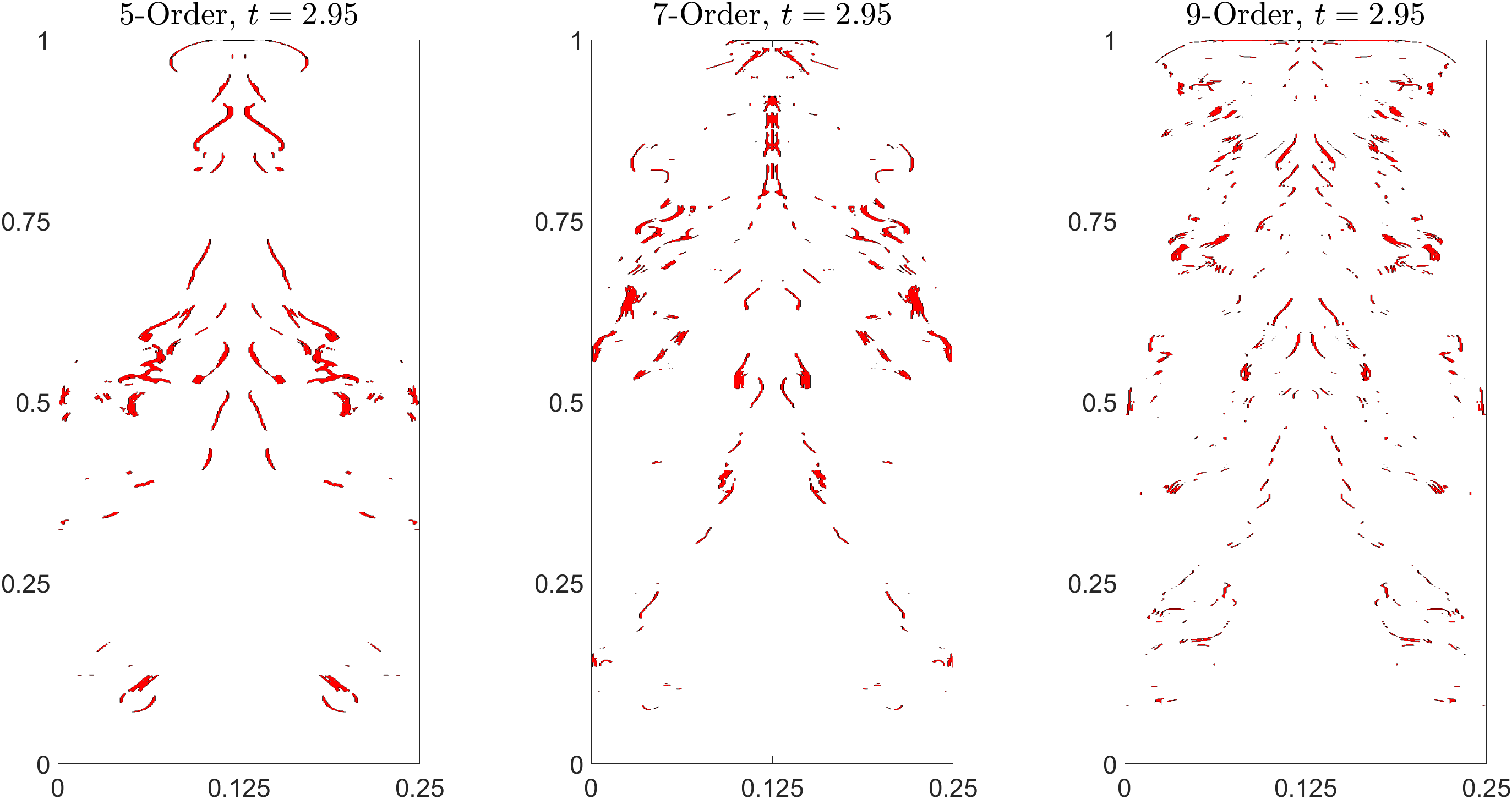}}
\caption{\sf Example 6: Same as in Figure \ref{fig10a}, but at the final time $t=2.95$.\label{fig10b}}
\end{figure}

We compute the numerical solution until the final time $t=2.95$ by the 5-Order (with the adaption constants $\texttt{C}_1=0.08$ and
$\texttt{C}_2=0.008$), 7-Order (with the adaption constants $\texttt{C}_1=0.1$ and $\texttt{C}_2=0.007$), and 9-Order (with the adaption
constants $\texttt{C}_1=0.12$ and $\texttt{C}_2=0.009$) schemes on a uniform mesh with $\dx=\dy=1/1024$ in the computational domain
$[0,0.25]\times[0,1]$. The obtained numerical results are presented in the top rows of Figures \ref{fig10a} and \ref{fig10b} at times
$t=1.95$ and $2.95$, respectively. The ``rough'' areas detected by the studied schemes at the same times are plotted in  the bottom rows of
Figures \ref{fig10a} and \ref{fig10b}. As one can see, the 7-Order and 9-Order schemes resolve more small solution structures: this
indicates that these schemes contain a smaller amount of numerical dissipation compared with their 5-Order counterpart.   
 
\begin{DA}
\paragraph{Funding.} The work of S. Chu was funded by the DFG--SPP 2183: Eigenschaftsgeregelte Umformprozesse with the Project(s)
HE5386/19-2,19-3 Entwicklung eines flexiblen isothermen Reckschmiedeprozesses f\"ur die eigenschaftsgeregelte Herstellung von
Turbinenschaufeln aus Hochtemperaturwerkstoffen (424334423) and by the Deutsche Forschungsgemeinschaft (DFG, German Research
Foundation)--SPP 2410 Hyperbolic Balance Laws in Fluid Mechanics: Complexity, Scales, Randomness (CoScaRa) within the Project(s) HE5386/27-1
(Zuf\"allige kompressible Euler Gleichungen: Numerik und ihre Analysis, 525853336). The work of V. A. Kolotilov and V. V. Ostapenko was
partially supported by the Russian Science Foundation (project number 26-11-00230). The work of A. Kurganov was supported in part by NSFC
grant W2431004.

\paragraph{Conflicts of interest.} On behalf of all authors, the corresponding author states that there is no conflict of interest.

\paragraph{Data and software availability.} The data that support the findings of this study and FORTRAN codes developed by the authors and
used to obtain all of the presented numerical results are available from the corresponding author upon reasonable request.
\end{DA}

\bibliographystyle{siamnodash}
\bibliography{Chu-Kurganov}
\end{document}